\documentclass[11pt]{amsart}
\usepackage{a4wide}
\usepackage[utf8]{inputenc}
\usepackage[T1]{fontenc}
\usepackage{amsmath,amssymb,amsthm}
\usepackage{enumitem}
\usepackage{xcolor}
\usepackage{tikz}
\usepackage{hyperref}

\setlist[enumerate]{label=(\roman*),ref=(\roman*)}

\theoremstyle{plain}
\newtheorem{theorem}{Theorem}[section]
\newtheorem{lemma}[theorem]{Lemma}
\newtheorem{proposition}[theorem]{Proposition}
\newtheorem{corollary}[theorem]{Corollary}
\theoremstyle{definition}
\newtheorem{definition}[theorem]{Definition}
\newtheorem{assumption}[theorem]{Assumption}
\newtheorem{example}[theorem]{Example}
\theoremstyle{remark}
\newtheorem{remark}[theorem]{Remark}

\numberwithin{equation}{section}

\usepackage{calc}
\usepackage{ragged2e}
\usepackage{setspace}
\newcommand{\FundingLogos}{%
  \raisebox{0pt}{\includegraphics[height=1.5cm]{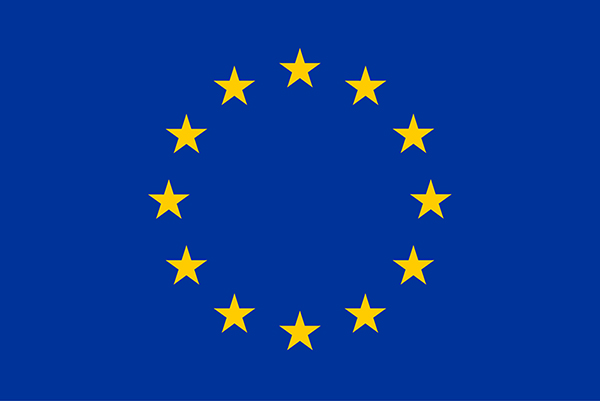}}%
  \hspace{1em}%
  \raisebox{0pt}{\includegraphics[height=1.5cm]{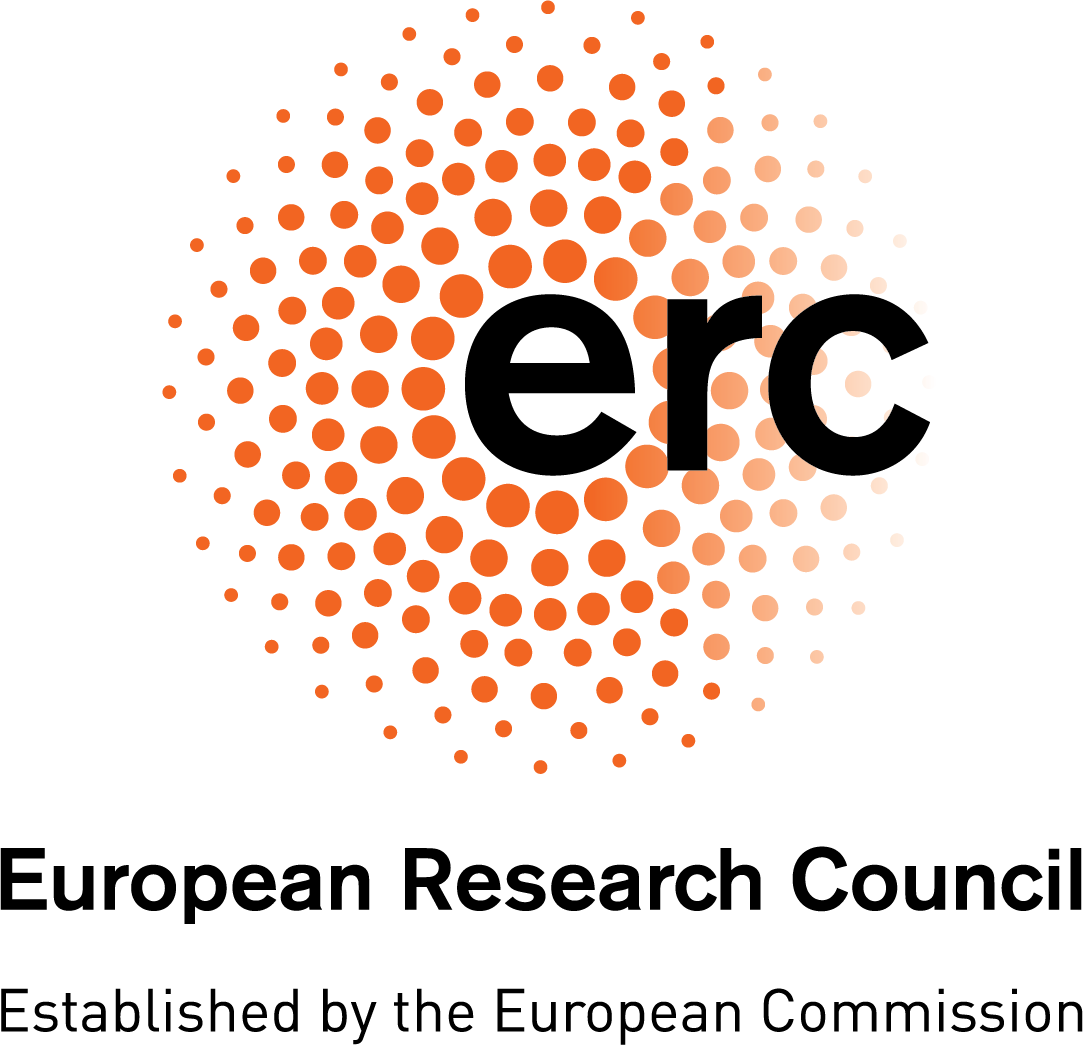}}%
}

\usepackage{tcolorbox}
\tcbset{colback=white,boxrule=0.5pt,arc=2pt,left=4pt,right=4pt,top=4pt,bottom=4pt}

\begin{document}

\setlength{\parindent}{0pt}

\title{On the continuity of martingale optimal transport in higher dimensions}

\author{Francesco Mattesini}
\address[Francesco Mattesini]{Department of Mathematics, CIT School, Technical University of Munich, \& Munich Center for Machine Learning (MCML), Germany}
\email{francesco.mattesini@tum.de}

\author{Johannes Wiesel}
\address[Johannes Wiesel]{Department of Mathematics, University of Copenhagen, Denmark}
\email{wiesel@math.ku.dk}

\date{\today}

\begin{abstract}
We investigate qualitative and quantitative continuity of martingale optimal transport under perturbations of both marginals in dimensions $d\geq2$. We assume that each pair of marginals admits a reference martingale coupling whose density relative to the product of its marginals is bounded above and away from zero, uniformly along the sequence. If the limiting target measure is not carried by a hyperplane, we prove continuity under $W_1$-convergence for continuous costs of at most linear growth, assuming only finite first moments.
For Lipschitz costs, a uniform $p$-moment bound on the targets,
$p>1$, yields an explicit stability estimate with exponent
$(1-1/p)/2$; uniformly bounded target supports give exponent $1/2$.
The proof constructs a correction of couplings with small martingale
bias while preserving both marginals. Truncated moment corrections
allow the target measure to have unbounded support. The construction also
shows that every martingale coupling of the limiting marginals can be approximated
in $W_1$ by martingale couplings of the approximating marginals.
We also discuss counterexamples showing that neither the uniform density condition
nor non-degeneracy of the limiting target measure can be omitted in general.
\end{abstract}

\maketitle

\section{Introduction}

Let $\mu$ and $\nu$ be two probability measures on $\mathbb{R}^d$ with finite first moment. A coupling $\pi(dx,dy)=\mu(dx)\pi_x(dy)$ of $\mu$ and $\nu$ is a martingale coupling if 
\[
\int y\,\pi_x(dy)=x \qquad \text{for } \mu\text{-a.e. } x.
\]
We denote the set of martingale couplings by $\Pi_M(\mu,\nu)$. 
By Strassen's theorem \cite[Theorem~8]{Str65} the set $\Pi_M(\mu,\nu)$ is nonempty if and only if $\mu$ and $\nu$ are in convex order, denoted by $\mu\preceq_c\nu$, which we assume throughout.
For a continuous cost $c:\mathbb{R}^d\times\mathbb{R}^d\to\mathbb{R}$ satisfying $|c(x,y)|\le C(1+|x|+|y|)$, the martingale optimal transport (MOT) problem is defined as
\begin{equation}\label{eq:MOTProb}
V(\mu,\nu):=\inf_{\pi\in\Pi_M(\mu,\nu)}\int c\,d\pi.
\end{equation}
In this paper we investigate continuity of $V$ with respect to the marginals $\mu,\nu$, i.e.~we ask if the functional $(\mu,\nu)\mapsto V(\mu,\nu)$ is continuous in the Wasserstein metric. In dimension $d=1$, such continuity holds, see \cite{BP22,Wie23}. In dimension $d\geq 2$ continuity fails in general, see \cite{BJ22}. 
The reason is that the set $\Pi_M(\mu_n,\nu_n)$ need not approximate $\Pi_M(\mu,\nu)$. Indeed, limits of martingale couplings in $\Pi_M(\mu_n,\nu_n)$ are martingale couplings in $\Pi_M(\mu,\nu)$, see Lemma~\ref{lem:compact}, but the converse may fail if the pairs $(\mu_n,\nu_n)$ are degenerate in a sense we specify below (cf.~Assumptions \ref{ass:1} and \ref{ass:2}). In the example of \cite{BJ22}, $\Pi_M(\mu_n,\nu_n)$ consists of a single coupling for every $n$, while $\Pi_M(\mu,\nu)$ contains many martingale couplings, some of them with strictly smaller cost. 
This is due to the fact that, for every $n$, the martingale constraint forces the mass to move along prescribed directions, and this constraint disappears in the limit.
The aim of this article is to identify a quantitative condition which prevents such hidden constraints, and under which $V$ is continuous with an explicit rate.

\subsection{Main result}
To exclude situations as in \cite{BJ22}, where every point can send its mass only to a few prescribed points of the target, we require the existence of one martingale coupling which spreads the mass of every point over the whole target, with weights comparable to those of the target.

\begin{assumption}\label{ass:1}
Let $\kappa\in(0,1)$. There exists $\bar\pi\in\Pi_M(\mu,\nu)$ with
\begin{equation}\label{ass:A}\tag{$A_\kappa$}
\kappa\leq\frac{d\bar\pi}{d(\mu\otimes\nu)}\leq\kappa^{-1}\qquad \mu\otimes\nu\text{-a.e.}
\end{equation}
\end{assumption}

By Fubini's theorem, \eqref{ass:A} means that the kernels of $\bar\pi$ satisfy
\[
\kappa\,\nu\leq\bar\pi_x\leq\kappa^{-1}\nu\qquad\text{for } \mu\text{-a.e. } x.
\]
We note that Assumption \ref{ass:1} asks for the existence of one coupling $\bar\pi$ satisfying \eqref{ass:A}. It is not a condition required for the whole set $\Pi_M(\mu,\nu)$. In particular the optimizer of $V(\mu,\nu)$ is typically singular with respect to $\mu\otimes\nu$. Assuming existence of martingale measures equivalent to $\mu\otimes\nu$ is classical in the literature on entropic martingale optimal transport, see e.g.~\cite{nutz2024martingale, chen2026convergence} for related variants.  As we will see below, \eqref{ass:A} alone is not sufficient for continuity of $V$, and it has to be complemented by a non-degeneracy condition on the target $\nu$.

\medskip

Let us now explain why this assumption is natural for the stability problem. In order to bound $V(\mu_n,\nu_n)$ from above, we construct specific martingale couplings $\pi_n$ of $(\mu_n,\nu_n)$, whose optimal transport costs $\int c\,d\pi_n$ are close to $V(\mu,\nu)$. This can be achieved by the shadow construction of \cite{EN22}, i.e.~by moving an optimizer of the limit problem onto the new marginals. Such a perturbation produces a coupling with the correct marginals and almost the correct cost. However, it is only an approximate martingale coupling. The question is therefore whether a coupling with a small martingale bias can be corrected into an exact martingale coupling of the same marginals at a small cost. The coupling $\bar\pi$ of \eqref{ass:A} is exactly the object which makes such a correction possible. Since $\bar\pi_x\geq\kappa\nu$, one may add to it signed perturbations of size smaller than $\kappa\nu$ without losing positivity. Since $\bar\pi_x\leq\kappa^{-1}\nu$, the errors created on the second marginal by such perturbations remain controlled. In other words, \eqref{ass:A} quantifies the room available to move mass while respecting the martingale constraint. A similar tilting construction was used in \cite[Ass. 2.1, Lemma 3.2]{nutz2024martingale} to show boundedness of dual optimizers, see also \cite[Ass. 2.3, Prop. 2.8]{chen2026convergence}.

\medskip

Assumption \eqref{ass:A} also implies that the pair $(\mu,\nu)$ is irreducible in the sense of \cite{DMT19}, see Remark~\ref{rem:irreducible}. It is however a strictly stronger requirement, since it forces the support of $\mu$ to stay at a quantitative distance from the relative boundary of the component. In particular $\text{supp}(\mu)$ has to be bounded, see Lemma~\ref{lem:support}. We think of \eqref{ass:A} as a quantitative version of irreducibility. For stability of $V$ it is required for the sequence $(\mu_n, \nu_n)$, with a constant $\kappa$ independent of $n$, and it then holds for the limit pair as well, see Lemma~\ref{lem:Alimit}. 
We now state a class of examples, where Assumption \ref{ass:1} is naturally satisfied.

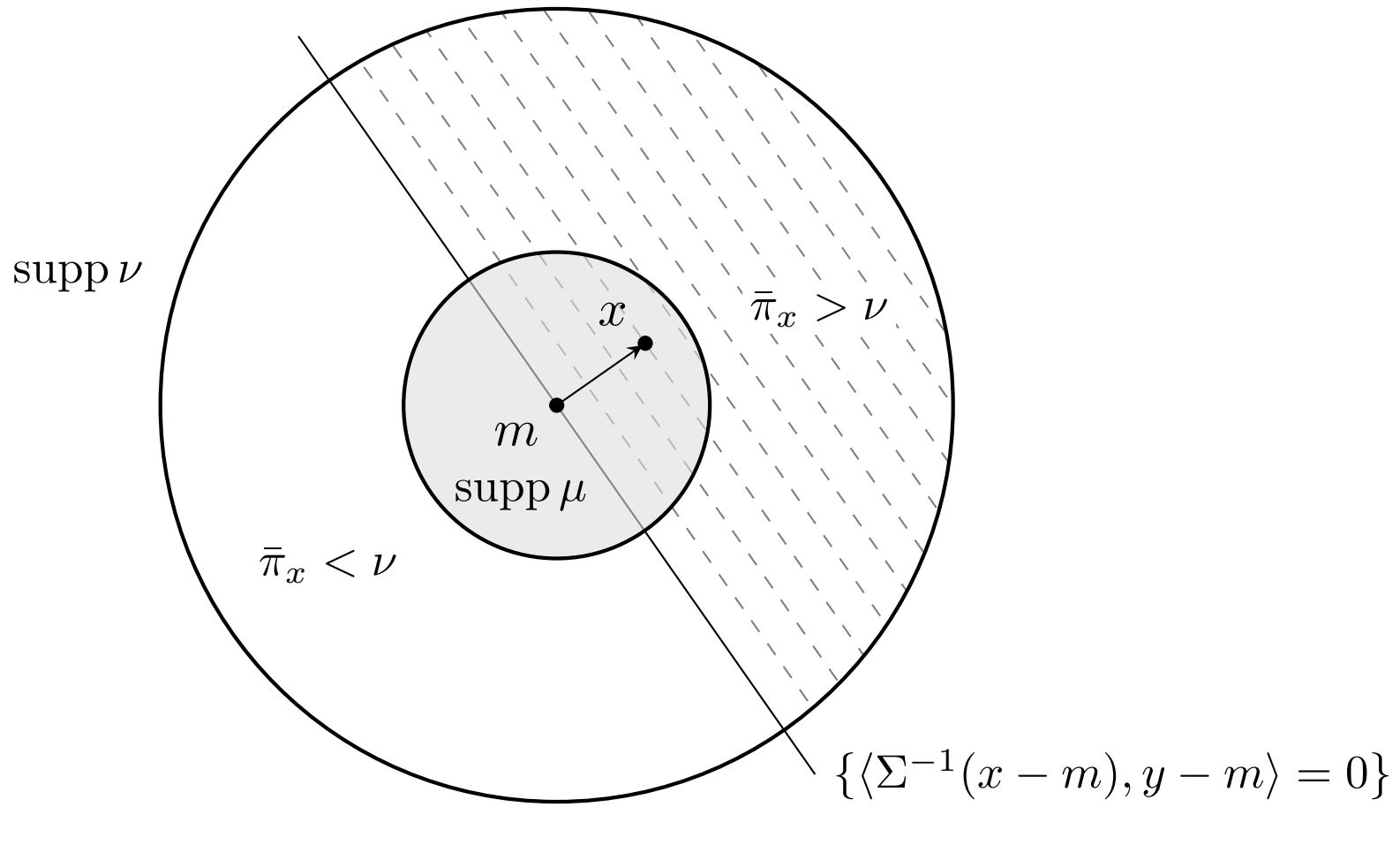
\begin{figure}[ht]
\centering
\begin{tikzpicture}[scale=1.35]
  \begin{scope}[rotate=35]
    \begin{scope}
      \clip (0,0) circle (2.2);
      \clip (0,-2.3) rectangle (2.3,2.3);
      \foreach \t in {0.2,0.4,...,2.2}
        \draw[dashed, gray] (\t,-2.3) -- (\t,2.3);
    \end{scope}
    \draw[thin] (0,-2.5) -- (0,2.5);
    \node[right, font=\small] at (0,-2.5) {$\{\langle \Sigma^{-1}(x-m),y-m\rangle=0\}$};
  \end{scope}
  \draw[thick] (0,0) circle (2.2);
  \draw[thick, fill=gray!25, fill opacity=0.6] (0,0) circle (0.85);
  \draw[->, >=stealth] (0,0) -- (35:0.58);
  \fill (0,0) circle (1.2pt) node[below left] {$m$};
  \fill (35:0.6) circle (1.2pt) node[above left] {$x$};
  \node[font=\small] at (-0.2,-0.5) {$\operatorname{supp}\mu$};
  \node[font=\small] at (165:2.75) {$\operatorname{supp}\nu$};
  \node[font=\small, fill=white, inner sep=1.5pt] at (20:1.55) {$\bar\pi_x>\nu$};
  \node[font=\small, fill=white, inner sep=1.5pt] at (215:1.55) {$\bar\pi_x<\nu$};
\end{tikzpicture}
\caption{The construction of Example~\ref{rem:examples}, with $\Sigma$ a multiple of the identity. The kernel $\bar\pi_x$ puts more mass than $\nu$ on the dashed half-space and less on the other half.}
\label{fig:tilt}
\end{figure}

\begin{example}\label{rem:examples}
Let $\nu$ have mean $m$, invertible covariance matrix $\Sigma$ and support in $B(m,R)$, where here and in the following we denote by $B(x,R)$ the closed ball of radius $R$ centered at $x$, and let $\mu$ have mean $m$ with
\[
R\,\big|\Sigma^{-1}(x-m)\big|\leq\theta<1\qquad\text{for } x\in\operatorname{supp}\mu.
\]
Then the kernel
\[
\bar\pi_x(dy):=\big(1+\langle\Sigma^{-1}(x-m),y-m\rangle\big)\,\nu(dy)
\]
has mean 
\begin{align*}
\int y\,\bar\pi_x(dy)&=m+\int(y-m)\langle \Sigma^{-1}(x-m),y-m\rangle\,\nu(dy)+m\,\Big\langle \Sigma^{-1}(x-m),\int(y-m)\,\nu(dy)\Big\rangle\\
&=m+\int(y-m)(y-m)^\top\,\nu(dy)\,\Sigma^{-1}(x-m) \\
& = m+\Sigma\Sigma^{-1}(x-m)=x.
\end{align*}
Since $\mu$ has mean $m$, integrating the kernels against $\mu$ gives
\[
\int\bar\pi_x\,\mu(dx)=\Big(1+\Big\langle\Sigma^{-1}\Big(\int x\,\mu(dx)-m\Big),y-m\Big\rangle\Big)\nu(dy)=\nu(dy),
\]
so that $\bar\pi(dx,dy):=\mu(dx)\,\bar\pi_x(dy)$ has second marginal $\nu$. Moreover, $|\langle\Sigma^{-1}(x-m),y-m\rangle|\leq\theta$ for $x\in\operatorname{supp}\mu$ and $y\in\operatorname{supp}\nu$, so that the density of $\bar\pi_x$ with respect to $\nu$ lies in $[1-\theta,1+\theta]$. Hence \eqref{ass:A} holds with $\kappa=1-\theta$. The construction is illustrated in Figure~\ref{fig:tilt}. Every kernel $\bar\pi_x$ is a tilt of $\nu$ which puts more mass than $\nu$ on the half-space pointing towards $x$, and less on the opposite one, in such a way that the barycenter moves from $m$ to $x$. An analogous construction, replacing the linear function $y$ by its truncation $\psi_T(y)$ defined in~\eqref{eq:phiT} gives examples with unbounded $\nu$, for instance a Gaussian $\nu$ and a uniform $\mu$ on a small ball.
\end{example}

We denote the $p$-th moment of $\nu$ by $M_p(\nu):=\int|y|^p\,d\nu$, and for sequences $\mu_n\preceq_c\nu_n$ converging to $\mu\preceq_c\nu$ we define the Wasserstein difference\footnote{Here and throughout we denote by $W_1$ the Wasserstein distance of order $1$, see \eqref{eq:Wass}.}
\[
\Delta_n:=W_1(\mu_n,\mu)+W_1(\nu_n,\nu).
\]
As we see below, only first moments are needed for continuity of $V$. Higher moments of $\nu_n$ quantify the tails and enter the rate, while no additional moment assumption on $\mu_n$ is needed since $\mu_n\preceq_c \nu_n$.
Before we state our main result, we introduce a second assumption. 

\begin{assumption}\label{ass:2}
The target measure $\nu$ is not carried by a hyperplane.
\end{assumption}
We call $\nu$ satisfying Assumption \ref{ass:2} non-degenerate.
\begin{theorem}\label{thm:main}
Let $\mu_n\preceq_c\nu_n$ be probability measures on $\mathbb{R}^d$ with finite first moment, converging in $W_1$ to $\mu\preceq_c\nu$. Assume that there is $\kappa\in(0,1)$ such that every pair $(\mu_n,\nu_n)$ satisfies \eqref{ass:A}, and that $\nu$ satisfies Assumption \ref{ass:2}.
\begin{enumerate}
\item\label{it:main-cont} If $c$ is continuous and there exists a constant $C > 0$\footnote{Here and throughout we denote by $C$ a positive constant, the value of which may change from line to line. $C$ may only depend on the dimension, the cost $c$, $\kappa$ and $\mu,\nu$ but not on $n$. For the quantitative estimates, $C$ may additionally depend on $p$ and $\sup_n M_p(\nu_n)$, or on a common support radius in the bounded-support case.} such that $|c(x,y)|\leq C(1+|x|+|y|)$, then $V(\mu_n,\nu_n)\to V(\mu,\nu)$.
\item\label{it:main-rate} If there exists $L>0$ such that $c$ is $L$-Lipschitz for the distance $|x-x'|+|y-y'|$, and $\sup_n M_p(\nu_n)<\infty$ for some $p\in(1,\infty)$, then there exists $n_0 > 0$ such that for every $n \ge n_0$
\[
\big|V(\mu_n,\nu_n)-V(\mu,\nu)\big|\leq L\,\Delta_n+C\,\Delta_n^{\gamma},\qquad \gamma:=\frac{1-1/p}{2},
\]
with $C$ independent of $n$. If the supports of all $\nu_n$ lie in a fixed bounded set, the same holds with $\gamma=\frac12$.
\end{enumerate}
\end{theorem}

Assumption \ref{ass:2} has a simple geometric interpretation: the support of the target measure should not become thin in any direction as $n\to\infty$. To see why this is important, recall that the martingale condition has to hold for every direction. In a direction in which $\nu_n$ keeps a fixed spread, this condition is flexible. A small error in the barycenter can be corrected by slightly reweighting $\nu_n$ along that direction.
In a direction where the support of $\nu_n$ becomes thin, the condition becomes rigid. There is almost no room to move mass in that direction, so balancing the barycenter now imposes relations between the remaining coordinates of the coupling. These constraints disappear in the limit, where the thin directions are lost. The set of martingale couplings then jumps, and the value \eqref{eq:MOTProb} can drop. In the example of \cite{BJ22} the martingale constraint forces transport along prescribed lines. 

\medskip

 Assumptions \ref{ass:1} and \ref{ass:2} play complementary roles. Assumption \ref{ass:1} ensures that every kernel of the reference coupling $\bar\pi$ explores the whole support of $\nu_n$, while Assumption \ref{ass:2} ensures that $\nu$ leaves room in every direction. Assumption \ref{ass:2} is placed on $\nu$ and not on $\mu$ for two reasons: the conditional laws we aim to produce are reweightings of $\nu$, and it is the spread of $\nu$ that provides room for such a procedure. Moreover, assuming non-degeneracy of the support of $\nu$ is weaker than assuming non-degeneracy of $\mu$, when $\mu,\nu$ have finite second moments. Indeed, every martingale coupling satisfies
\[
\Sigma(\nu)=\Sigma(\mu)+\int\operatorname{Cov}(\pi_x)\,\mu(dx).
\]
In consequence, non-degeneracy of $\mu$ implies non-degeneracy of $\nu$.

\medskip

Below we give a number of examples, showing that neither Assumption \ref{ass:1} nor Assumption \ref{ass:2} alone is sufficient in general for Theorem \ref{thm:main}. In Example~\ref{ex:parabola} we assume that $\nu_n$ lives on a thin band around a parabola, so that it approximately fixes the conditional second moment of every kernel. This happens although there is $\kappa\in(0,1)$ such that every pair satisfies Assumption \ref{ass:1}. Assumption \ref{ass:2} excludes such a mechanism. It guarantees that, for $n$ large, $\nu_n$ has spread in every direction, uniformly in $n$, see Lemma~\ref{lem:MT}.

\medskip
Conversely, Assumption \ref{ass:2} alone is not sufficient either. To see this one can employ the variation of the counterexample of \cite{BJ22} given in \cite[Remark~3.2(ii)]{BJ22}, in which a small multiple of a measure with full support is added to both marginals. The limit target then has full support, while the same argument as in \cite[Proposition~1.2(ii)]{BJ22} shows that $V$ is discontinuous. By Theorem~\ref{thm:main}, these pairs cannot satisfy~\eqref{ass:A} with a constant $\kappa$ independent of $n$.

\medskip
Finally, the non-degeneracy of $\nu$ is not implied by \eqref{ass:A}. The latter is invariant under affine bijections applied to both marginals, while flattening maps make the covariance degenerate. For a single pair $(\mu,\nu)$ imposing Assumption \ref{ass:2} is harmless, since one may work in the affine hull of $\operatorname{supp}\nu$, which also carries $\mu$. However it poses a constraint on the sequences $(\mu_n,\nu_n)$, excluding those with a loss of dimension in the limit.
\subsection{Main steps in the proof of Theorem~\ref{thm:main}}\label{sec:proof.main}

In the following we describe the main steps in the proof of Theorem~\ref{thm:main}. For the detailed proofs we refer to Sections~\ref{sec:shadow}, \ref{sec:repair} and~\ref{sec:proof}.

\medskip

Our starting point is the well-known observation that optimal martingale couplings of $(\mu_n,\nu_n)$ are relatively compact, and that their limit points are martingale couplings of $(\mu,\nu)$. This gives the lower bound
\[
\liminf_{n\to\infty}V(\mu_n,\nu_n)\geq V(\mu,\nu)
\]
for every sequence, without any assumption, see Corollary~\ref{cor:lsc}. The main difficulty thus lies in proving the upper bound. We will achieve this by constructing good martingale couplings of $(\mu_n,\nu_n)$.

\medskip

For this we first relax the martingale constraint via penalization. For a coupling $\pi$ of $(\mu,\nu)$ we measure how far it is from being a martingale coupling by the $L^1(\mu)$ norm of its \emph{martingale bias}
\[
b_\pi(x):=\int(y-x)\,\pi_x(dy),
\]
and we consider the penalized optimal transport problem
\begin{equation}\label{eq:PR}
P_R(\mu,\nu):=\inf_{\pi\in\Pi(\mu,\nu)}\Big\{\int c\,d\pi+R\,\|b_\pi\|_{L^1(\mu)}\Big\}.
\end{equation}
Its constraint set consists of all couplings of $\mu$ and $\nu$, and therefore behaves well under perturbations of the marginals. We highlight two properties: first $P_R\leq V$, since martingale couplings have no bias. Second, $P_R$ is Lipschitz w.\,r.\,t.~the marginals if $c$ is $L$-Lipschitz, with Lipschitz constant $(R+L)$, see Proposition~\ref{prop:PRstab}.

\medskip

Let us comment on the choice of penalty. A natural alternative would have been to add to the cost a Lagrangian term
\[
\int\langle h(x),y-x\rangle\,\pi(dx,dy)=\int\langle h,b_\pi\rangle\,d\mu
\]
for a fixed multiplier $h:\mathbb{R}^d\to\mathbb{R}^d$. Such a term is not a penalty. It vanishes on martingale couplings, but it can be negative, and it does not see any martingale bias orthogonal to $h$. To recover $V$ from it, one has to take the supremum over all multipliers and exchange it with the infimum over couplings, see e.g.~the seminal work \cite{BJ16}. This exchange requires a minimax theorem, hence some compactness, and the optimal multipliers are in general unbounded. For our purposes, the term $R\,\|b_\pi\|_{L^1(\mu)}$ is somewhat more natural. It is obtained by taking the supremum over bounded multipliers before taking the infimum over couplings, i.e.~
\[
\sup_{|h|\leq R}\int\langle h(x),y-x\rangle\,\pi(dx,dy)=R\,\|b_\pi\|_{L^1(\mu)},
\]
see Remark~\ref{rem:lagrange}. In this form no exchange of infimum and supremum is needed. The penalty is nonnegative, it vanishes exactly on martingale couplings, and it requires no regularity of the multiplier. A fixed $h$ would enter the stability estimate through its Lipschitz constant, while $\|b_\pi\|_{L^1(\mu)}$ is controlled under perturbations of the marginals by the Wasserstein distance alone. Finally, the $L^1$-norm of the martingale bias is the natural choice here. A perturbation of the marginals of size $\Delta$ in $W_1$ creates a bias of size $\Delta$ in $L^1$, see Lemma~\ref{lem:shadow}, but in general there is no control in $L^\infty$.

\medskip

The next step is to prove stability for the penalized problem $P_R$. 
As mentioned above, our main tool for this is the shadow of a coupling introduced in \cite[Definition~3.1]{EN22}, which moves a coupling onto new marginals coordinate by coordinate. The new ingredient introduced here is an estimate on the martingale bias of the shadow, see Lemma~\ref{lem:shadow}. It states that moving a martingale coupling by an amount $\Delta$ in the Wasserstein sense creates a bias of size at most $\Delta$ in $L^1$.

\medskip

Since the constant in the stability estimate of $P_R$ grows with $R$, the stability of $V=\lim_R P_R$ requires an additional input, which constitutes the main novelty of our proof. We show that under the assumptions of Theorem~\ref{thm:main}, $V(\mu_n,\nu_n)-P_R(\mu_n,\nu_n)$ converges to zero as $R\to\infty$, uniformly in $n\ge n_0$; see Corollary~\ref{cor:uniform}. This is deduced from the fact that a coupling with a small martingale bias can be repaired into a martingale coupling at a small cost, see Proposition~\ref{prop:repair}. The repair proceeds in three moves, all of which use the coupling $\bar\pi$ of \eqref{ass:A} as a reservoir. First, at the few points where the bias is large, we replace the kernel of the coupling by the kernel of $\bar\pi$, which has no bias. Second, this replacement changes the second marginal by a signed measure of zero mass, which we absorb by mixing in a small amount of $\bar\pi$ and subtracting a product correction. Third, the remaining bias is now uniformly small, and we cancel it exactly by a tilt of $\nu$ along the direction of the bias, see Lemma~\ref{lem:smalldefect}. When $\nu$ has unbounded support we use a truncated tilt. Whether or not its support is bounded, this is where the support of $\nu$ and its non-degeneracy enter, see Lemma~\ref{lem:MT}.

\medskip

Combining the three steps and assuming that $c$ is $L$-Lipschitz, we obtain
\[
V(\mu_n,\nu_n)\leq P_R(\mu_n,\nu_n)+e(R)\leq P_R(\mu,\nu)+(L+R)\Delta_n+e(R)\leq V(\mu,\nu)+(L+R)\Delta_n+e(R),
\]
where $e(R)$ is the uniform error in the comparison between $V$ and $P_R$. The term $R\Delta_n$ is the price of a large penalty, and $e(R)$ the price of a small one. Balancing the two terms gives the rate of Theorem~\ref{thm:main}, see Section~\ref{sec:proof}. We then show in Section~\ref{sec:example} that the non-degeneracy of $\nu$ cannot be dropped.

\subsection{Related literature}

Martingale optimal transport initially emerged from ideas in model-independent pricing and hedging. Foundational discrete-time and continuous-time formulations can be found in \cite{BHLP13} and \cite{GHT14}. Subsequent work developed its variational and geometric structure, including the left-curtain coupling of \cite{BJ16}, complete quasi-sure duality on the line \cite{BNT17}, and, in higher dimensions, the decomposition of martingale transports into irreducible components \cite{GKL19, OS17, DMT19}. \cite{GO19} developed convergent discretization schemes for multistep, multidimensional MOT using a relaxation of
the martingale constraint. Regarding stability of $V$, \cite{BP22} and \cite{Wie23} established general continuity results in dimension one, see also \cite{BJMP22} for approximation of arbitrary martingale couplings in adapted Wasserstein distance. For $d>1$, \cite{BJ22} showed that both continuity of $V$ and stability of optimizers can fail, providing the main motivation for this paper.

\medskip
The related work \cite{CMS} establishes constraint qualification and dual attainment for transport problems with moment constraints. Very recently, \cite{CXZ} establishes quantitative $W_2$-stability for Wasserstein-relaxed entropic martingale optimal transport and convergence to classical MOT, with quantitative approximation bounds under assumptions similar to our Assumptions \ref{ass:1} and \ref{ass:2}. Their formulation accommodates marginals outside convex order and their construction uses moment perturbations related to ours. Our contribution was obtained independently and concerns continuity of the unregularized MOT value with both marginals prescribed.

\medskip
Entropic relaxation of the MOT problem has also received considerable attention in the recent decade.
\cite{DM18} studied entropic approximation and numerical
algorithms in several dimensions, while
\cite{HL19} introduced martingale Schr\"odinger bridges for the
calibration of stochastic volatility models. 
In a similar spirit, \cite{Kay25} proposed martingale couplings minimizing the relative entropy with respect to the product of their marginals as a method to calibrate multi-asset models to option prices.
For prescribed one-dimensional marginals,
\cite{NWZ23a,NWZ23b,nutz2024martingale} established the existence of
Schr\"odinger potentials and dual attainment, and \cite{chen2026convergence} proved exponential
convergence of a martingale Sinkhorn algorithm. More recently,
\cite{KRW26} obtained multidimensional exponential Sinkhorn convergence under compactness and geometric
assumptions. Our continuity results for unregularized MOT under
perturbations of both marginals complement these developments.

\section{Preliminaries}\label{sec:prelim}
\subsection{Notation}
Throughout the paper we denote by $\mathcal{P}_1(\mathbb{R}^k)$ the probability measures on $\mathbb{R}^k$ with finite first moment. For $\mu,\nu\in\mathcal{P}_1(\mathbb{R}^d)$ we denote by $\Pi(\mu,\nu)$ the set of couplings of $\mu$ and $\nu$. Every $\pi\in\Pi(\mu,\nu)$ disintegrates as
\[
\pi(dx,dy)=\mu(dx)\,\pi_x(dy),
\]
and we denote by $\pi_x$  the conditional distribution of $\pi$ given $x$. On $\mathbb{R}^d\times\mathbb{R}^d$ we use the distance $|x-x'|+|y-y'|$, and we identify $\mathcal{P}_1(\mathbb{R}^d\times\mathbb{R}^d)$ with the probability measures $\pi$ such that $\int(|x|+|y|)\,d\pi<\infty$. Every coupling of two measures in $\mathcal{P}_1(\mathbb{R}^d)$ belongs to this set.

\medskip

For two measures $\mu,\tilde\mu\in\mathcal{P}_1(\mathbb{R}^d)$ we denote by
\begin{equation}\label{eq:Wass}
W_1(\mu,\tilde\mu):=\inf_{\pi\in\Pi(\mu,\tilde\mu)}\int_{\mathbb{R}^d\times\mathbb{R}^d}|x-x'|\,\pi(dx,dx')
\end{equation}
the Wasserstein distance of order one. For two couplings $\pi,\tilde\pi\in\mathcal{P}_1(\mathbb{R}^d\times\mathbb{R}^d)$ we denote by
\[
W_1(\pi,\tilde\pi):=\inf_{\Gamma\in\Pi(\pi,\tilde\pi)}\int_{(\mathbb{R}^d\times\mathbb{R}^d)^2}\big(|x-x'|+|y-y'|\big)\,\Gamma(dx,dy,dx',dy')
\]
the Wasserstein distance of order one with respect to the distance $|x-x'|+|y-y'|$, where $\Pi(\pi,\tilde\pi)$ denotes the couplings of $\pi$ and $\tilde\pi$ as probability measures on $\mathbb{R}^d\times\mathbb{R}^d$.

\medskip

Both infima are attained, see \cite[Theorem~4.1]{Vil09}. Note that if $\pi\in\Pi(\mu,\nu)$ and $\tilde\pi\in\Pi(\tilde\mu,\tilde\nu)$, then projecting a coupling of $\pi$ and $\tilde\pi$ onto the two pairs of first and second coordinates gives
\[
W_1(\mu,\tilde\mu)+W_1(\nu,\tilde\nu)\leq W_1(\pi,\tilde\pi).
\]
Recall that two measures are in convex order, $\mu\preceq_c\nu$, if and only if  for every convex function $f:\mathbb{R}^d\to\mathbb{R}$, 
\[
\int f\,d\mu\leq\int f\,d\nu.
\]
Moreover, by Strassen's theorem \cite[Theorem~8]{Str65}, the set of martingale couplings $\Pi_M(\mu,\nu)$ is not empty if and only if $\mu\preceq_c\nu$.

\medskip

For $\pi\in\Pi(\mu,\nu)$ we call
\begin{equation}\label{eq:defect}
b_\pi(x):=\int(y-x)\,\pi_x(dy)
\end{equation}
the martingale bias of $\pi$, so that $\pi\in\Pi_M(\mu,\nu)$ if and only if $b_\pi=0$ for $\mu$-a.e.\ $x$. We measure its size by $\|b_\pi\|_{L^1(\mu)}$ and $\|b_\pi\|_{L^\infty(\mu)}$.
For a finite signed measure $\zeta$ on $\mathbb{R}^k$ we denote by
\[
\|\zeta\|_{TV}:=|\zeta|(\mathbb{R}^k)
\]
the total variation norm. We will repeatedly use the following characterization of convergence in $W_1$, see \cite[Definition~6.8 and Theorem~6.9]{Vil09}: a sequence $(\mu_n)$ in $\mathcal{P}_1(\mathbb{R}^k)$ converges to $\mu$ in $W_1$ if and only if it converges weakly and its first moments converge, i.e.~
\[
\int|z|\,\mu_n(dz)\to\int|z|\,\mu(dz).
\]
In this case $\int f\,d\mu_n\to\int f\,d\mu$ for every continuous $f$ satisfying $|f(z)|\leq C(1+|z|)$. In particular, if $|c(x,y)|\leq C(1+|x|+|y|)$, the map $\pi\mapsto\int c\,d\pi$ is continuous with respect to $W_1$.

\subsection{Basic facts on martingale couplings}
We start with an elementary identity, which is the only property of the convex order used in the construction described in Section \ref{sec:proof.main}. We emphasize that it holds for every coupling, not only for martingale couplings.

\begin{lemma}\label{lem:zeromean}
If $\mu\preceq_c\nu$, then for every $\pi\in\Pi(\mu,\nu)$
\begin{equation}\label{eq:zeromean}
\int b_\pi\,d\mu=0.
\end{equation}
\end{lemma}
\begin{proof}
The functions $x\mapsto\pm x_j$ are convex, hence by convex order we have $\int x\,d\mu=\int y\,d\nu$. Moreover,
\[
\int b_\pi\,d\mu=\int\!\!\int y\,\pi_x(dy)\,\mu(dx)-\int x\,d\mu=\int y\,d\nu-\int x\,d\mu=0.\qedhere
\]
\end{proof}

Next we show that the martingale constraint is closed under $W_1$-convergence of the marginals. This yields the existence of optimizers and the lower semicontinuity of $V$.

\begin{lemma}\label{lem:compact}
Let $\mu_n\to\mu$ and $\nu_n\to\nu$ in $W_1$, and $\pi_n\in\Pi(\mu_n,\nu_n)$. Then
\begin{enumerate}
\item\label{it:compact-rel} $(\pi_n)$ is relatively compact in $W_1$, and every limit point lies in $\Pi(\mu,\nu)$,
\item\label{it:compact-mart} if $\pi_n\in\Pi_M(\mu_n,\nu_n)$ for all $n$, every limit point lies in $\Pi_M(\mu,\nu)$.
\end{enumerate}
\end{lemma}

\begin{proof}
For (i), the families $(\mu_n)$ and $(\nu_n)$ are tight, so $(\pi_n)$ is tight, see \cite[Lemma~4.4]{Vil09}. Let $\pi$ be a weak limit point along a subsequence. Its marginals are $\mu$ and $\nu$, since the projections are continuous. Moreover
\[
\int(|x|+|y|)\,d\pi_n=\int|x|\,d\mu_n+\int|y|\,d\nu_n\to\int|x|\,d\mu+\int|y|\,d\nu=\int(|x|+|y|)\,d\pi,
\]
so that the convergence holds in $W_1$.

\medskip

For (ii), fix $\varphi\in C_b(\mathbb{R}^d)$, $j\leq d$, and set $f(x,y):=\varphi(x)(y_j-x_j)$. This function is continuous with $|f|\leq\|\varphi\|_\infty(|x|+|y|)$, and
\[
\int f\,d\pi_n=\int\varphi\,(b_{\pi_n})_j\,d\mu_n=0.
\]
By $W_1$-convergence we obtain $\int\varphi\,(b_\pi)_j\,d\mu=\int f\,d\pi=0$. The finite signed measure $(b_\pi)_j\,d\mu$ thus vanishes on $C_b(\mathbb{R}^d)$, hence it is zero, and $b_\pi=0$ for $\mu$-a.e.\ $x$.
\end{proof}

\begin{corollary}\label{cor:lsc}
Let $c$ be continuous with $|c(x,y)|\leq C(1+|x|+|y|)$. The infimum defining $V(\mu,\nu)$ is attained, and
\[
\liminf_{n\to\infty}V(\mu_n,\nu_n)\geq V(\mu,\nu)
\]
whenever $\mu_n\to\mu$ and $\nu_n\to\nu$ in $W_1$.
\end{corollary}

\begin{proof}
We first show that the infimum is attained. Applying Lemma~\ref{lem:compact} to the constant sequences $\mu_n=\mu$ and $\nu_n=\nu$, we see that $\Pi_M(\mu,\nu)$ is compact in $W_1$. Since $|c(x,y)|\leq C(1+|x|+|y|)$, the map $\pi\mapsto\int c\,d\pi$ is continuous with respect to $W_1$. A continuous function on a compact set attains its minimum, so there is an optimizer for $V(\mu,\nu)$.

\medskip

We now turn to the lower bound. For every $n$, let $\pi_n\in\Pi_M(\mu_n,\nu_n)$ be an optimizer for $V(\mu_n,\nu_n)$, which exists by the first part of the proof. We pass to a subsequence, which we do not relabel, along which $V(\mu_n,\nu_n)$ converges to $\liminf_n V(\mu_n,\nu_n)$. By Lemma~\ref{lem:compact}, a further subsequence of $(\pi_n)$ converges in $W_1$ to some $\pi\in\Pi_M(\mu,\nu)$. Using once more the continuity of $\pi\mapsto\int c\,d\pi$, we conclude that
\[
\liminf_{n\to\infty}V(\mu_n,\nu_n)=\lim_{n\to\infty}\int c\,d\pi_n=\int c\,d\pi\geq V(\mu,\nu),
\]
where the last inequality holds because $\pi$ is a martingale coupling of $\mu$ and $\nu$.
\end{proof}

\subsection{\texorpdfstring{Consequences of~\eqref{ass:A}}{Consequences of (A kappa)}}
For $\nu\in\mathcal{P}_1(\mathbb{R}^d)$ we write
\begin{equation}\label{eq:mean}
m(\nu):=\int y\,d\nu,\qquad
\widehat M_1(\nu):=\int|y-m(\nu)|\,d\nu.
\end{equation}
The first consequence of \eqref{ass:A} is that $\mu$ lives in a ball around the mean of $\nu$, whose radius is controlled by $\kappa$. Indeed, since every kernel is dominated by $\kappa^{-1}\nu$, it cannot put enough mass far away to balance a point which is far from the mean.

\begin{lemma}\label{lem:support}
If $(\mu,\nu)$ satisfies \eqref{ass:A}, then $|x-m(\nu)|\leq\kappa^{-1}\widehat M_1(\nu)$ for $\mu$-a.e.\ $x$.
\end{lemma}
\begin{proof}
Write $m:=m(\nu)$ and let $\bar\pi$ be the martingale coupling satisfying~\eqref{ass:A}. For $\mu$-a.e.\ $x$, the kernel $\bar\pi_x$ has barycenter $x$ and satisfies $\bar\pi_x\leq\kappa^{-1}\nu$. Fix such an $x$. We may assume w.l.o.g.~that 
$x\neq m$.

\medskip

Since $\bar\pi_x$ is a probability measure with barycenter $x$, we can write the distance of $x$ from $m$ as an integral against $\bar\pi_x$, namely
\[
|x-m|=\left\langle \frac{x-m}{|x-m|} ,x-m\right\rangle=\left\langle \frac{x-m}{|x-m|},\int(y-m)\,\bar\pi_x(dy)\right\rangle=\int\left\langle \frac{x-m}{|x-m|},y-m\right\rangle\,\bar\pi_x(dy).
\]
Finally, by the latter equality together with \eqref{ass:A} we have
\[
|x-m|= \int\left\langle \frac{x-m}{|x-m|},y-m\right\rangle\,\bar\pi_x(dy)\stackrel{\eqref{ass:A}}{\leq}
\kappa^{-1}\int|y-m|\,\nu(dy) \stackrel{\eqref{eq:mean}}{=} \kappa^{-1}\widehat M_1(\nu).
\]
\end{proof}

The lower bound in \eqref{ass:A} has a geometric consequence as well. It forces every kernel to reach the whole support of $\nu$, which identifies the irreducible components.

\begin{remark}[Relation with irreducibility]\label{rem:irreducible}
Assumption~\eqref{ass:A} implies that the pair $(\mu,\nu)$ is irreducible, with a single component in the sense of \cite{DMT19}. To make this precise, for a set $A\subset\mathbb{R}^d$ we denote by $\operatorname{conv}A$ its convex hull, by $\operatorname{cl}A$ its closure, and by $\operatorname{ri}A$ its relative interior, i.e., its interior relative to the affine hull of $A$, see \cite[Section~6]{Roc70}. Then the irreducible convex paving $I$ of \cite{DMT19} satisfies
\[
I(x)=\operatorname{ri}\operatorname{cl}\operatorname{conv}\operatorname{supp}\nu\qquad\text{for }\mu\text{-a.e. } x.
\]
To see this, let $\bar\pi$ be the coupling satisfying~\eqref{ass:A}. On the one hand, \cite[Theorem~2.1(ii)]{DMT19} gives $I(x)\subset\operatorname{cl}\operatorname{conv}\operatorname{supp}\nu$. On the other hand, \cite[Equation~(2.2)]{DMT19} gives $\bar\pi_x(\operatorname{cl}I(x))=1$. Since $\bar\pi_x\geq\kappa\nu$, the measure $\nu$ is absolutely continuous with respect to $\bar\pi_x$, so that $\nu(\operatorname{cl}I(x))=1$ as well. Hence $\operatorname{supp}\nu\subset\operatorname{cl}I(x)$, and since $I(x)$ is convex, $\operatorname{cl}\operatorname{conv}\operatorname{supp}\nu\subset\operatorname{cl}I(x)$. The two inclusions show that $I(x)$ and $\operatorname{conv}\operatorname{supp}\nu$ have the same closure. As $I(x)$ is relatively open, \cite[Corollary~6.3.1]{Roc70} gives the claim.

\medskip

However, Assumption~\eqref{ass:A} is stronger than irreducibility as it also keeps $\mu$ away from the relative boundary of the component. Indeed, since $\bar\pi_x\geq\kappa\nu$, we may write
\[
\bar\pi_x=\kappa\,\nu+(1-\kappa)\,q_x,\qquad q_x:=\frac{\bar\pi_x-\kappa\nu}{1-\kappa},
\]
where $q_x$ is a probability measure dominated by a multiple of $\nu$. 
Integrating $y$ with respect to both sides, and using that $\bar\pi_x$ has mean $x$ and $\nu$ has mean $m(\nu)$, we obtain
\[
x-m(\nu)=(1-\kappa)\Big(\int y\,q_x(dy)-m(\nu)\Big)\in(1-\kappa)\big(\operatorname{cl}\operatorname{conv}\operatorname{supp}\nu-m(\nu)\big)
\]
for $\mu$-a.e.\ $x$. Since the right-hand side is a closed set, this gives
\[
\operatorname{supp}\mu\subset m(\nu)+(1-\kappa)\big(\operatorname{cl}\operatorname{conv}\operatorname{supp}\nu-m(\nu)\big).
\]
This inclusion fails for irreducible pairs such as $\mu$ the uniform measure on $[-1,1]$ and $\nu=\frac12(\delta_{-1}+\delta_1)$, for which the support of $\mu$ reaches the boundary of $\operatorname{cl}\operatorname{conv}\operatorname{supp}\nu=[-1,1]$.
\end{remark}

The next lemma shows that \eqref{ass:A} with a fixed $\kappa$ is stable under limits of the marginals in the $W_1$-metric.
\begin{lemma}\label{lem:Alimit}
Let $\mu_n\to\mu$ and $\nu_n\to\nu$ in the Wasserstein metric $W_1$, and assume every $(\mu_n,\nu_n)$ satisfies \eqref{ass:A}. Then $(\mu,\nu)$ satisfies \eqref{ass:A}.
\end{lemma}

\begin{proof}
{\sc Step 1.} (Convergence of the products). The measures $\mu_n\otimes\nu_n$ converge weakly to $\mu\otimes\nu$. Indeed the family is tight, and every limit point $\pi$ of $(\mu_n\otimes \nu_n)$ satisfies
\[
\int f(x)g(y)\,d\pi=\int f\,d\mu\int g\,d\nu\qquad\text{for } f,g\in C_b(\mathbb{R}^d).
\]
These functions determine $\pi$ on the class of measurable rectangles, which is closed under intersections, hence $\pi=\mu\otimes\nu$.

\medskip

{\sc Step 2.} (Passage to the limit). Let $\bar\pi_n\in\Pi_M(\mu_n,\nu_n)$ be the couplings satisfying~\eqref{ass:A}. By Lemma~\ref{lem:compact}, a subsequence converges in $W_1$ to some $\bar\pi\in\Pi_M(\mu,\nu)$. The density bounds in~\eqref{ass:A} give, for every $f\in C_b(\mathbb{R}^d\times\mathbb{R}^d)$ with $f\geq0$,
\[
\kappa\int f\,d(\mu_n\otimes\nu_n)\leq\int f\,d\bar\pi_n\leq\kappa^{-1}\int f\,d(\mu_n\otimes\nu_n).
\]
By Step 1 and the weak convergence of $\bar\pi_n$, we may let $n\to\infty$ along the subsequence, and obtain
\[
\kappa\int f\,d(\mu\otimes\nu)\leq\int f\,d\bar\pi\leq\kappa^{-1}\int f\,d(\mu\otimes\nu).
\]
Finally, two finite Borel measures on a metric space which are ordered on nonnegative bounded continuous functions are ordered as measures, by approximating indicators of closed sets from above. Hence $\kappa\,\mu\otimes\nu\leq\bar\pi\leq\kappa^{-1}\mu\otimes\nu$, i.e., $\bar\pi$ has a density with respect to $\mu\otimes\nu$ with values in $[\kappa,\kappa^{-1}]$, and $(\mu,\nu)$ satisfies~\eqref{ass:A}.
\end{proof}
\subsection{Truncated second moments and tails}

The repair map described in Section \ref{sec:proof.main} corrects a martingale bias by tilting $\nu$ along the direction of the bias. However, when the support of $\nu$ is unbounded, the linear tilt $\langle v,y-m\rangle\,\nu(dy)$ need not be dominated by a multiple of $\nu$. 
To overcome this difficulty we introduce a truncated tilt. We fix 
a continuous nonincreasing $\chi:[0,\infty)\to[0,1]$ with $\chi=1$ on $[0,1]$ and $\chi=0$ on $[2,\infty)$. For $T>0$ and $\nu\in\mathcal{P}_1(\mathbb{R}^d)$ with mean $m$, let
\begin{equation}\label{eq:phiT}
\psi_T(y):=(y-m)\,\chi\bigg(\frac{|y-m|}T \bigg),\qquad \varphi_T:=\psi_T-\int\psi_T\,d\nu,
\end{equation}
and let
\begin{equation}\label{eq:MT}
\Sigma_T(\nu):=\int(y-m)(y-m)^\top\chi\bigg(\frac{|y-m|}T \bigg)\,\nu(dy),\qquad \lambda_T(\nu):=\lambda_{\min}\big(\Sigma_T(\nu)\big).
\end{equation}
Then, by construction
\begin{equation}\label{eq:est.phiT.MT}
\text{
$|\psi_T|\leq2T$, $|\varphi_T|\leq4T$, and $\Sigma_T(\nu)$ is symmetric and positive semidefinite.
}
\end{equation} 
The matrix $\Sigma_T(\nu)$ plays the role of the covariance of $\nu$, and its smallest eigenvalue measures the spread of $\nu$ in its thinnest direction.

\begin{lemma}\label{lem:MT}
Let $\nu\in\mathcal{P}_1(\mathbb{R}^d)$ and $T>0$, and let $\varphi_T$, $\Sigma_T(\nu)$ and $\lambda_T(\nu)$ be defined in~\eqref{eq:phiT} and~\eqref{eq:MT}. Then the following hold.
\begin{enumerate}
\item\label{it:MT-id} For every $v\in\mathbb{R}^d$,
\[
\int\langle v,\varphi_T\rangle\,d\nu=0\qquad\text{and}\qquad\int y\,\langle v,\varphi_T(y)\rangle\,\nu(dy)=\Sigma_T(\nu)v.
\]
\item\label{it:MT-pos} If $\nu$ is not carried by a hyperplane, then $\lambda_T(\nu)>0$ for $T$ large enough.
\item\label{it:MT-cont} If $\nu_n\to\nu$ in $W_1$, then $\Sigma_T(\nu_n)\to \Sigma_T(\nu)$ and $\lambda_T(\nu_n)\to\lambda_T(\nu)$ for every $T>0$.
\end{enumerate}
\end{lemma}

\begin{proof}
{\sc Proof of }\ref{it:MT-id}. By the definition~\eqref{eq:phiT}, $\varphi_T$ is obtained from $\psi_T$ by subtracting its mean, hence
\begin{equation}\label{eq:trunc.first.id}
\int\langle v,\varphi_T\rangle\,d\nu=\Big\langle v,\int\psi_T\,d\nu-\int\psi_T\,d\nu\Big\rangle=0.
\end{equation}
For the second identity we write $y=(y-m)+m$. 
Combining the latter with \eqref{eq:trunc.first.id} yields
\begin{align*}
\int y\,\langle v,\varphi_T(y)\rangle\,\nu(dy)&
\stackrel{\eqref{eq:phiT},\eqref{eq:trunc.first.id}}{=}\int(y-m)\langle v,\psi_T(y)\rangle\,\nu(dy)-\Big\langle v,\int\psi_T\,d\nu\Big\rangle\int(y-m)\,\nu(dy)\\
&\stackrel{\eqref{eq:mean}, \eqref{eq:phiT}}{=}\int(y-m)(y-m)^\top\chi\bigg(\frac{|y-m|}T\bigg)\,\nu(dy)\,v\stackrel{\eqref{eq:MT}}{=}\Sigma_T(\nu)v.
\end{align*}

\medskip

{\sc Proof of }\ref{it:MT-pos}. For a unit vector $e$ we consider the quadratic form of $\Sigma_T(\nu)$ in the direction $e$, and we define
\begin{align}\label{eq:help}
f_T(e):=\langle \Sigma_T(\nu)e,e\rangle\stackrel{\eqref{eq:MT}}{=}\int\langle e,y-m\rangle^2\chi\bigg(\frac{|y-m|}T\bigg)\,\nu(dy),
\end{align}
so that 
\[
\lambda_T(\nu)=\min_{|e|=1}f_T(e).
\]
By \eqref{eq:phiT} and \eqref{eq:est.phiT.MT}, the integrand in \eqref{eq:help} is bounded by $4T^2$ and continuous in $e$, so that $f_T$ is continuous by dominated convergence. Moreover, $f_T$ is nondecreasing in $T$, because $\chi$ is nonincreasing. 

\medskip
We assume by contradiction 
that $\lambda_T(\nu)=0$ for arbitrarily large $T$. By monotonicity, $\lambda_T(\nu)=0$ for every $T>0$. Choose $T_k\to\infty$ and unit vectors $e_k$ with $f_{T_k}(e_k)=0$. By compactness of the unit sphere we may assume w.l.o.g that $e_k\to e$ for some unit vector $e$. For fixed $T$ and $k$ large enough so that $T_k\geq T$, monotonicity gives $f_T(e_k)\leq f_{T_k}(e_k)=0$, and by continuity $f_T(e)=0$. Letting $T\to\infty$, monotone convergence gives
\[
\int\langle e,y-m\rangle^2\,\nu(dy)=0.
\]
Since the integrand is nonnegative, $\langle e,y-m\rangle=0$ for $\nu$-a.e.\ $y$, so that $\nu$ is carried by the hyperplane $\{y:\langle e,y-m\rangle=0\}$, a contradiction.

\medskip

{\sc Proof of }\ref{it:MT-cont}. Recall that convergence of $\nu_n\to\nu$ in $W_1$ implies weak convergence and convergence of the moments $m_n:=m(\nu_n)\to m$. Let $g(z):=zz^\top\chi(|z|/T)$. Then, $g(z)$ is continuous with compact support, hence bounded and uniformly continuous. By~\eqref{eq:MT} we may write
\[
\Sigma_T(\nu_n)=\int g(y-m_n)\,\nu_n(dy),\qquad\Sigma_T(\nu)=\int g(y-m)\,\nu(dy).
\]
Since $m_n\to m$, uniform continuity gives 
\[
\sup_y|g(y-m_n)-g(y-m)| \to 0.
\]
Since $y\mapsto g(y-m)$ is bounded and continuous, weak convergence gives 
\[
\int g(y-m)\,\nu_n(dy)\to\int g(y-m)\,\nu(dy).
\]
Putting the last two together yields
\[
\big|\Sigma_T(\nu_n)-\Sigma_T(\nu)\big|\leq\sup_y\big|g(y-m_n)-g(y-m)\big|+\Big|\int g(y-m)\,\nu_n(dy)-\int g(y-m)\,\nu(dy)\Big|\to0.
\]
Finally, recall that the smallest eigenvalue of a matrix is $1$-Lipschitz on symmetric matrices, that is 
$
|\lambda_{\min}(A)-\lambda_{\min}(B)|\leq\|A-B\|,
$
for all symmetric matrices $A, B$. Thus, $\lambda_T(\nu_n)\to\lambda_T(\nu)$.
\end{proof}

When estimating the repair error we need to control the first moments
of submeasures of $\nu$ with small mass. Such submeasures may
concentrate in the tails of $\nu$; the required control is provided
by the following modulus of uniform integrability:
\begin{equation}\label{eq:omega}
\omega_\nu(a):=\sup\Big\{\int|y|\,d\rho:\ \rho\leq\nu,\ \rho(\mathbb{R}^d)\leq a\Big\},\qquad a\geq0,
\end{equation}
which is nondecreasing in $a$. For every $t>0$ and every $\rho\leq\nu$ of mass at most $a$,
\[
\int|y|\,d\rho\leq\int_{\{|y|>t\}}|y|\,d\nu+ta,
\]
so that 
\begin{equation}\label{eq:buon.unif.int}
\omega_\nu(a)\leq\int_{\{|y|>t\}}|y|\,d\nu+ta.
\end{equation}
If $\nu_n\to\nu$ in $W_1$, the family $\{|y|\}$ is uniformly integrable under $(\nu_n)$, and \eqref{eq:buon.unif.int} shows that
\begin{equation}\label{eq:omegaUI}
\lim_{a\to0}\,\sup_n\omega_{\nu_n}(a)=0.
\end{equation}
If $M_p(\nu)<\infty$ with $p>1$, by Hölder's inequality we obtain
\begin{equation}\label{eq:omegaHolder}
\omega_\nu(a)\leq M_p(\nu)^{\frac 1 p }a^{1-\frac 1 p}.
\end{equation}
Finally, if $\operatorname{supp}\nu\subset B(0,R)$,
then $\omega_\nu(a)\leq R a$.

\section{Shadows of couplings and the penalized value}\label{sec:shadow}
We recall the notion of a shadow of a coupling from \cite[Definition~3.1]{EN22}: given a coupling $\pi$ of $(\mu,\nu)$ and new marginals $\tilde\mu$ and $\tilde\nu$, the shadow of $\pi$ is a coupling of $(\tilde\mu,\tilde\nu)$ obtained by moving the first marginal to $\tilde\mu$ and the second marginal to $\tilde\nu$, each along a $W_1$-optimal coupling, see Definition~\ref{def:shadow}.
It should not be confused with the shadow of a measure in the sense of \cite{BJ16}, which is not used in this paper.

\begin{definition}\label{def:shadow}
Let $\mu,\nu,\tilde\mu,\tilde\nu\in\mathcal{P}_1(\mathbb{R}^d)$. Let $\kappa_1\in\Pi(\mu,\tilde\mu)$ and $\kappa_2\in\Pi(\nu,\tilde\nu)$ be couplings attaining $W_1(\mu,\tilde\mu)$ and $W_1(\nu,\tilde\nu)$, and let $\kappa_1=\mu\otimes K_1$ and $\kappa_2=\nu\otimes K_2$ be their disintegrations. Given $\pi\in\Pi(\mu,\nu)$, its shadow $\tilde\pi\in\Pi(\tilde\mu,\tilde\nu)$ is defined as the second marginal of $\pi\otimes K\in\mathcal{P}\big((\mathbb{R}^d\times\mathbb{R}^d)^2\big)$, where the kernel $K$ is defined as $K(x,y):=K_1(x)\otimes K_2(y)$.
\end{definition}

The existence of optimal couplings is granted by \cite[Theorem~4.1]{Vil09}. However, $W_1$-optimal couplings need not be unique, so the shadow of $\pi$ may depend on the choice of $\kappa_1$ and $\kappa_2$. All the statements below hold for every such choice. We therefore fix one choice and speak of the shadow, as in \cite{EN22}.

\medskip

In the proofs we make use of the following probabilistic representation: let $(X,Y,X',Y')$ be random variables with joint law
\[
\pi\otimes K(dx,dy,dx',dy')=\pi(dx,dy)\,K_1(x,dx')\,K_2(y,dy').
\]
Then the shadow $\tilde\pi$ is the law of $(X',Y')$. The product form of $\pi\otimes K$ gives three facts, which we use repeatedly:
\begin{enumerate}[label=(S\arabic*),ref=(S\arabic*)]
\item\label{it:fact-Y'} Given $(X,Y,X')$, the law of $Y'$ is $K_2(Y,\cdot)$.
\item\label{it:fact-Y} Given $(X,X')$, the law of $Y$ is $\pi_X$. Indeed, the joint law of the triple $(X,Y,X')$ is $\mu(dx)\,\pi_x(dy)\,K_1(x,dx')$, in which $y$ and $x'$ are conditionally independent given $x$.
\item\label{it:fact-pairs} The pair $(X,X')$ has law $\mu(dx)K_1(x,dx')$ and the pair $(Y,Y')$ has law $\nu(dy)K_2(y,dy')$. Hence $\mathbb{E}|X-X'|=W_1(\mu,\tilde\mu)$ and $\mathbb{E}|Y-Y'|=W_1(\nu,\tilde\nu)$.
\end{enumerate}

The next lemma collects the properties of the shadow that we need. The first two items say that the shadow $\tilde\pi$ is a coupling of the new marginals which is close to $\pi$, both in $W_1$ and in terms of the 
{\color{blue}
transport
}
cost. The third item is the new ingredient. It shows that if $\pi$ is a martingale coupling, then its shadow has martingale bias of at most $\Delta$ in $L^1$. In general we cannot expect the shadow to be a martingale coupling, because $X$ and $Y$ are moved independently of each other.

\begin{lemma}\label{lem:shadow}
Let $\Delta:=W_1(\mu,\tilde\mu)+W_1(\nu,\tilde\nu)$. Then
\begin{enumerate}
\item\label{it:shadow-dist} $\tilde\pi\in\Pi(\tilde\mu,\tilde\nu)$ and $W_1(\tilde\pi,\pi)\leq\Delta$,
\item\label{it:shadow-cost} $\big|\int c\,d\tilde\pi-\int c\,d\pi\big|\leq L\Delta$ if $c$ is $L$-Lipschitz,
\item\label{it:shadow-defect} $\|b_{\tilde\pi}\|_{L^1(\tilde\mu)}\leq\|b_\pi\|_{L^1(\mu)}+\Delta$.
\end{enumerate}
\end{lemma}

\begin{proof}
Item~\ref{it:shadow-dist} is \cite[Lemma~3.2]{EN22} with $p=1$ and two marginals, where it is shown that in fact $W_1(\tilde\pi,\pi)=\Delta$. Item~\ref{it:shadow-cost} follows from \ref{it:shadow-dist}, since $c$ is $L$-Lipschitz for the distance $|x-x'|+|y-y'|$. The new statement is \ref{it:shadow-defect}.

\medskip

{\sc Proof of \ref{it:shadow-defect}.} Let $(X,Y,X',Y')$ be the random variables introduced above. Since $\tilde\pi$ is the law of $(X',Y')$, the kernel of $\tilde\pi$ at $X'$ is the conditional law of $Y'$ given $X'$, so that
\[
b_{\tilde\pi}(X')=\mathbb{E}\big[Y'-X'\,\big|\,X'\big]\qquad\text{almost surely.}
\]
We want to compare this quantity with the martingale bias $b_\pi(X)$ of $\pi$. To this end, we write $Y'-X'$ as the sum of increments. We first show that, after conditioning on $X'$, the middle increment gives exactly $\mathbb{E}[b_\pi(X)|\,X']$. Then we show that the two remaining increments are controlled by the displacements $W_1(\nu,\tilde\nu)$ and $W_1(\mu,\tilde\mu)$.

\medskip

We start by replacing $Y'$ by its conditional mean. Let $\hat Y(y):=\int y'\,K_2(y,dy')$. By the tower property and \ref{it:fact-Y'},
\[
b_{\tilde\pi}(X')=\mathbb{E}\Big[\mathbb{E}\big[Y'\,\big|\,X,Y,X'\big]-X'\,\Big|\,X'\Big]\stackrel{\ref{it:fact-Y'}}{=}\mathbb{E}\big[\hat Y(Y)-X'\,\big|\,X'\big].
\]
Splitting $\hat Y(Y)-X'=(\hat Y(Y)-Y)+(Y-X)+(X-X')$ and using the tower property once more for the middle term, we obtain
\[
b_{\tilde\pi}(X')=\mathbb{E}\big[\hat Y(Y)-Y\,\big|\,X'\big]+\mathbb{E}\Big[\mathbb{E}\big[Y-X\,\big|\,X,X'\big]\,\Big|\,X'\Big]+\mathbb{E}\big[X-X'\,\big|\,X'\big].
\]
By \ref{it:fact-Y}, the inner conditional expectation of the middle term is the bias of $\pi$, namely
\[
\mathbb{E}\big[Y-X\,\big|\,X,X'\big]\stackrel{\ref{it:fact-Y}}{=}\int(y-X)\,\pi_X(dy)\stackrel{\eqref{eq:defect}}{=}b_\pi(X).
\]
Hence
\[
b_{\tilde\pi}(X')=\mathbb{E}\big[\hat Y(Y)-Y\,\big|\,X'\big]+\mathbb{E}\big[b_\pi(X)\,\big|\,X'\big]+\mathbb{E}\big[X-X'\,\big|\,X'\big].
\]

\medskip

It remains to estimate the three terms. Taking norms and expectations, and using Jensen's inequality for conditional expectations, we get
\begin{equation}\label{eq:defect.formula.proof}
\begin{split}
\|b_{\tilde\pi}\|_{L^1(\tilde\mu)}=\mathbb{E}|b_{\tilde\pi}(X')|&\leq\mathbb{E}|\hat Y(Y)-Y|+\mathbb{E}|b_\pi(X)|+\mathbb{E}|X-X'|\\
&\stackrel{\ref{it:fact-pairs}}{=}\mathbb{E}|\hat Y(Y)-Y|+\|b_\pi\|_{L^1(\mu)}+W_1(\mu,\tilde\mu),
\end{split}
\end{equation}
where we also used that $X$ has law $\mu$. For the first term, $\hat Y(y)-y=\int(y'-y)\,K_2(y,dy')$ gives $|\hat Y(y)-y|\leq\int|y'-y|\,K_2(y,dy')$, so that
\[
\mathbb{E}|\hat Y(Y)-Y|\leq\int\!\!\int|y'-y|\,K_2(y,dy')\,\nu(dy)=\mathbb{E}|Y'-Y|\stackrel{\ref{it:fact-pairs}}{=}W_1(\nu,\tilde\nu).
\]
Inserting this into \eqref{eq:defect.formula.proof} yields 
\[
\|b_{\tilde\pi}\|_{L^1(\tilde\mu)}\leq\|b_\pi\|_{L^1(\mu)}+W_1(\mu,\tilde\mu)+W_1(\nu,\tilde\nu),
\]
which is \ref{it:shadow-defect}.
\end{proof}

Using Lemma~\ref{lem:shadow}, we now deduce the stability of the penalized value~\eqref{eq:PR}.

\begin{proposition}\label{prop:PRstab}
Let $R\geq0$, and let $\mu_n\preceq_c\nu_n$ converge in $W_1$ to $\mu\preceq_c\nu$.
\begin{enumerate}
\item\label{it:PR-cont} If $c$ is continuous with $|c(x,y)|\leq C(1+|x|+|y|)$, then $\limsup_n P_R(\mu_n,\nu_n)\leq P_R(\mu,\nu)$.
\item\label{it:PR-lip} If $c$ is $L$-Lipschitz, then for all pairs
\[
\big|P_R(\tilde\mu,\tilde\nu)-P_R(\mu,\nu)\big|\leq(L+R)\big(W_1(\mu,\tilde\mu)+W_1(\nu,\tilde\nu)\big).
\]
\end{enumerate}
\end{proposition}

\begin{proof}
{\sc Proof of }\ref{it:PR-cont}. Fix an arbitrary $\eta>0$ and $\pi\in\Pi(\mu,\nu)$ with 
\[
\int c\,d\pi+R\|b_\pi\|_{L^1(\mu)}\leq P_R(\mu,\nu)+\eta,
\]
and let $\pi_n$ be its shadow on $(\mu_n,\nu_n)$. Since $W_1(\mu_n,\mu)\to0$ and $W_1(\nu_n,\nu)\to0$, Lemma~\ref{lem:shadow}\ref{it:shadow-dist} gives $W_1(\pi_n,\pi)\to0$. As $|c(x,y)|\leq C(1+|x|+|y|)$, the map $\pi\mapsto\int c\,d\pi$ is continuous with respect to $W_1$, so that $\int c\,d\pi_n\to\int c\,d\pi$. Moreover, Lemma~\ref{lem:shadow}\ref{it:shadow-defect} gives 
\[
\|b_{\pi_n}\|_{L^1(\mu_n)}\leq\|b_\pi\|_{L^1(\mu)}+\Delta_n,
\]
where $\Delta_n:=W_1(\mu_n,\mu)+W_1(\nu_n,\nu)\to0$. Hence
\[
\limsup_{n\to\infty}P_R(\mu_n,\nu_n)\leq\limsup_{n\to\infty}\Big(\int c\,d\pi_n+R\|b_{\pi_n}\|_{L^1(\mu_n)}\Big)\leq\int c\,d\pi+R\|b_\pi\|_{L^1(\mu)}\leq P_R(\mu,\nu)+\eta.
\]
Since $\eta>0$ is arbitrary, the claim follows.

\medskip

{\sc Proof of }\ref{it:PR-lip}. We argue as in the proof of~\ref{it:PR-cont}. Fix $\eta>0$ and $\pi\in\Pi(\mu,\nu)$ as above, let $\tilde\pi$ be its shadow on $(\tilde\mu,\tilde\nu)$, and set $\Delta:=W_1(\mu,\tilde\mu)+W_1(\nu,\tilde\nu)$. By Lemma~\ref{lem:shadow}\ref{it:shadow-cost} and~\ref{it:shadow-defect},
\[
P_R(\tilde\mu,\tilde\nu)\leq\int c\,d\tilde\pi+R\|b_{\tilde\pi}\|_{L^1(\tilde\mu)}\leq\int c\,d\pi+L\Delta+R\big(\|b_\pi\|_{L^1(\mu)}+\Delta\big)\leq P_R(\mu,\nu)+\eta+(L+R)\Delta.
\]
Letting $\eta\to0$, we get $P_R(\tilde\mu,\tilde\nu)-P_R(\mu,\nu)\leq(L+R)\Delta$. Exchanging the roles of the two pairs gives the same bound for $P_R(\mu,\nu)-P_R(\tilde\mu,\tilde\nu)$ and hence \ref{it:PR-lip}.
\end{proof}

\begin{remark}[Lagrange multipliers]\label{rem:lagrange}
For continuous $h:\mathbb{R}^d\to\mathbb{R}^d$ with $\sup|h|\leq R$ one has
\[
\int\langle h(x),y-x\rangle\,d\pi=\int\langle h,b_\pi\rangle\,d\mu\leq R\,\|b_\pi\|_{L^1(\mu)},
\]
with equality in the supremum over such $h$ by Lusin's theorem, see \cite[Theorem~2.24]{Rud87}. For compactly supported marginals, Sion's theorem, see \cite[Corollary~3.3]{Sio58}, then identifies $P_R(\mu,\nu)$ with
\[
\sup_{\sup|h|\leq R}\ \inf_{\pi\in\Pi(\mu,\nu)}\int\big(c(x,y)+\langle h(x),y-x\rangle\big)\,d\pi.
\]
Thus $P_R$ is the Lagrangian value with multipliers restricted to a ball, and Corollary~\ref{cor:uniform} says that the multipliers may be taken bounded, uniformly along the sequence, up to an error $e(R)$. 
\end{remark}

\section{The repair map}\label{sec:repair}

Proposition~\ref{prop:PRstab} shows that $P_R$ is stable under perturbations of the marginals, but with a constant of order $R$, so it does not by itself yield the stability of $V$. In this section we show that $V$ and $P_R$ can be compared uniformly along the sequence $(\mu_n,\nu_n)$, by repairing couplings with a small martingale bias. For this purpose, a single pair $\mu\preceq_c\nu$ in $\mathcal{P}_1(\mathbb{R}^d)$ satisfying Assumption \ref{ass:1} is fixed, with coupling $\bar\pi$ satisfying \eqref{ass:A}. We write $m:=m(\nu)$, fix $T>0$ with $\lambda:=\lambda_T(\nu)>0$, and write $\Sigma_T:=\Sigma_T(\nu)$, where $\Sigma_T(\nu)$ and $\lambda_T(\nu)$ are defined in~\eqref{eq:MT} and $\varphi_T$ is defined in~\eqref{eq:phiT}. We also set
\begin{equation}\label{eq:rAomega}
r:=\operatorname*{ess\,sup}_\mu|x|,\qquad
A:=1+r+\int|y|\,\nu(dy),\qquad
\omega:=\omega_\nu.
\end{equation}
Here $\omega_\nu$ is defined in~\eqref{eq:omega}, and $r<\infty$ by Lemma~\ref{lem:support}. Every $\pi\in\Pi(\mu,\nu)$ satisfies $\int(1+|x|+|y|)\,\pi(dx,dy)\leq A$.

\subsection{Tilting the target measure}

Let $v:\mathbb{R}^d\to\mathbb{R}^d$ be bounded and Borel. By Lemma~\ref{lem:MT}\ref{it:MT-id}, the signed kernel
\begin{equation}\label{eq:tilt}
\tau_x(dy):=\langle v(x),\varphi_T(y)\rangle\,\nu(dy)
\end{equation}
has zero mass and first moment $\Sigma_Tv(x)$. Integrated against $\mu$ it gives $\langle\int v\,d\mu,\varphi_T(y)\rangle\,\nu(dy)$, which vanishes when $\int v\,d\mu=0$. So $\tau$ shifts the conditional mean at $x$ by $\Sigma_Tv(x)$ without changing the marginals. The direction of the tilt is given by $v(x)$.

\begin{figure}[h!]
\centering
\begin{tikzpicture}[scale=1.35]
  \begin{scope}[rotate=30]
    \begin{scope}
      \clip (0,-2.3) rectangle (2.3,2.3);
      \clip (0,0) circle (2);
      \foreach \t in {0.15,0.3,...,2.2}
        \draw[dotted, gray] (\t,-2.3) -- (\t,2.3);
      \fill[white] (0,0) circle (1);
    \end{scope}
    \begin{scope}
      \clip (0,-2.3) rectangle (2.3,2.3);
      \clip (0,0) circle (1);
      \foreach \t in {0.15,0.3,...,2.2}
        \draw[dashed, gray] (\t,-2.3) -- (\t,2.3);
    \end{scope}
    \draw[thin] (0,-2.4) -- (0,2.4);
    \node[right, font=\small] at (0,-2.4) {$\{\langle v(x),y-m\rangle=0\}$};
  \end{scope}
  \draw[thick] (0,0) circle (1);
  \draw[thick, dashed] (0,0) circle (2);
  \draw[->, >=stealth, very thick] (0,0) -- (30:0.8);
  \node[font=\small, fill=white, inner sep=1pt] at (60:0.62) {$v(x)$};
  \fill (0,0) circle (1.2pt) node[below left] {$m$};
  \node[font=\small] at (250:1.17) {$T$};
  \node[font=\small] at (250:2.2) {$2T$};
  \node[font=\small, fill=white, inner sep=1.5pt] at (10:1.5) {$\tau_x>0$};
  \node[font=\small, fill=white, inner sep=1.5pt] at (200:1.45) {$\tau_x<0$};
  \node[font=\small] at (140:2.65) {$\tau_x=0$};
\end{tikzpicture}
\caption{The tilt~\eqref{eq:tilt} for $\nu$ symmetric around $m$. It is positive on the hatched half-space, linear inside the ball of radius $T$ (dashed), cut off by $\chi$ between radius $T$ and $2T$ (dotted), and zero outside the ball of radius $2T$.}
\label{fig:truncatedtilt}
\end{figure}
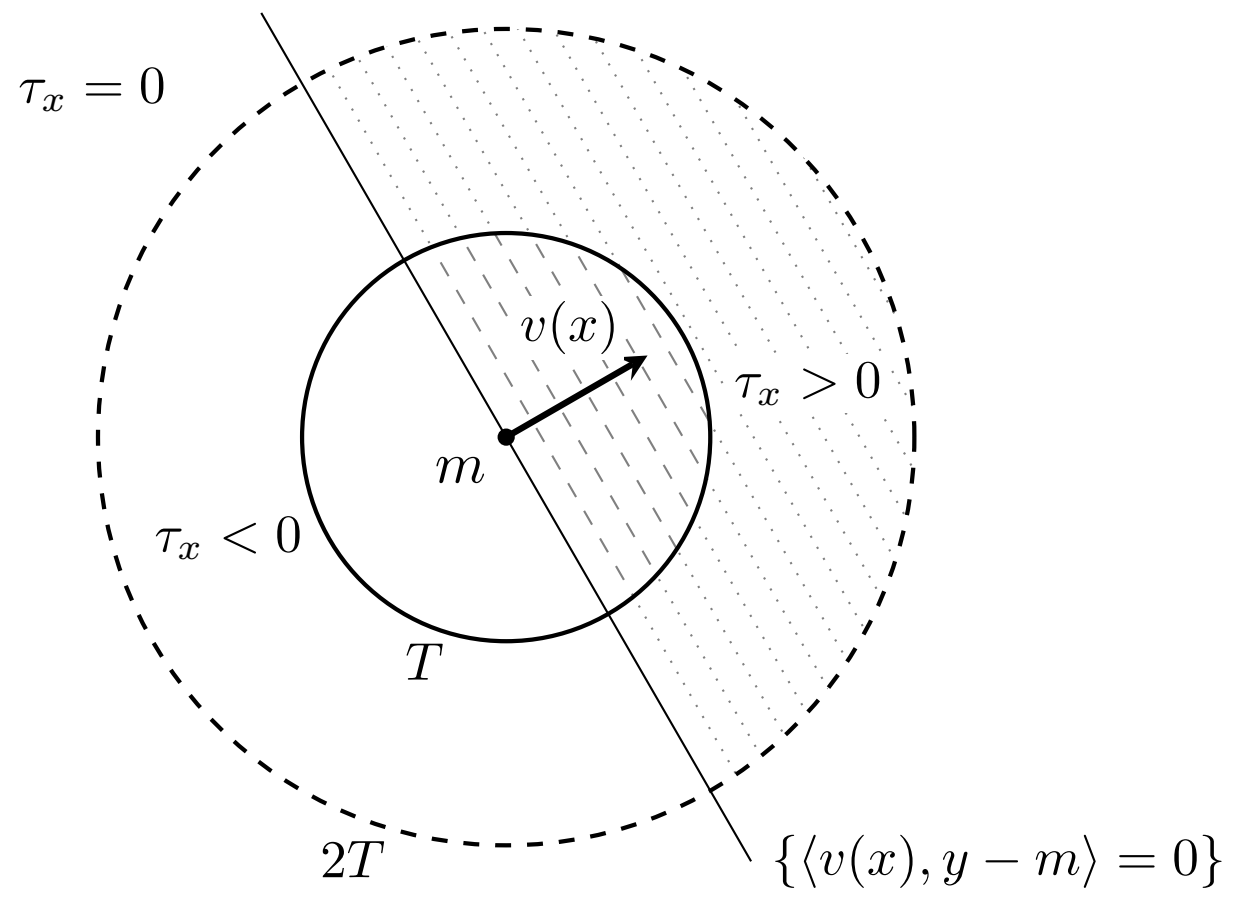

The density of the tilt is
\[
\frac{d\tau_x}{d\nu}(y)
=
\chi\Big(\frac{|y-m|}{T}\Big)\langle v(x),y-m\rangle
-
\Big\langle v(x),\int\psi_T\,d\nu\Big\rangle.
\]
If $\nu$ is symmetric around $m$, the centering term vanishes.
In this case the density is linear inside $B(m,T)$, is cut off
between radii $T$ and $2T$, and vanishes outside $B(m,2T)$;
see Figure~\ref{fig:truncatedtilt}.
For general $\nu$, centering introduces a constant offset.
Outside $B(m,2T)$ the density equals
$-\langle v(x),\int\psi_T\,d\nu\rangle$ and therefore need not
vanish. 

\medskip

The size of the tilt is controlled pointwise by
\begin{equation}\label{eq:taubound}
|\tau_x|\leq4T\,|v(x)|\,\nu.
\end{equation}
This domination is crucial for the construction below, as it bounds the negative part of the signed tilt by
a multiple of $\nu$. It is also the reason for the truncation by $\chi$. If we used the linear function $y-m$ in place of $\varphi_T$, the kernel need not be dominated by a multiple of $\nu$ when $\operatorname{supp}\nu$ is unbounded.

\subsection{Repair of a uniformly small martingale bias}

To repair the bias $b_\pi$ of a coupling $\pi$, we subtract the tilt~\eqref{eq:tilt} in the direction $v=\Sigma_T^{-1}b_\pi$, which shifts the conditional mean at $x$ by $b_\pi(x)$. 
Subtracting the tilt may create negative mass, since $\tau_x$ is a signed measure. To avoid this, we mix $\pi$ with a small amount of $\bar\pi$. Since $\bar\pi_x\geq\kappa\nu$ and $|\tau_x|\leq4T|v(x)|\,\nu$ by~\eqref{eq:taubound}, this mixing absorbs the negative part of $-\tau_x$ as soon as $v$ is small enough. We formalize this idea in the following lemma.

\begin{lemma}\label{lem:smalldefect}
Let $\bar\pi\in\Pi_M(\mu,\nu)$ be a martingale coupling satisfying \eqref{ass:A}, so that $\kappa\,\nu\leq\bar\pi_x\leq\kappa^{-1}\nu$ for $\mu$-a.e.\ $x$. Let $\pi\in\Pi(\mu,\nu)$ satisfy $\|b_\pi\|_{L^\infty(\mu)}\leq\kappa\lambda/(4T)$, and let $\tau$ be the tilt defined in~\eqref{eq:tilt} with $v:=\Sigma_T^{-1}b_\pi$. Set
\begin{equation}\label{eq:s}
s:=\frac{4T}{\kappa\lambda}\,\|b_\pi\|_{L^\infty(\mu)}\in[0,1]
\end{equation}
and, for $\mu$-a.e.\ $x$,
\[
\tilde\pi_x:=(1-s)\,\pi_x+s\,\bar\pi_x-(1-s)\,\tau_x,
\]
and $\tilde\pi(dx,dy):=\mu(dx)\,\tilde\pi_x(dy)$. Then $\tilde\pi\in\Pi_M(\mu,\nu)$ and
\[
\int(1+|x|+|y|)\,|\tilde\pi-\pi|(dx,dy)\leq A\Big(2s+\frac{4T}{\lambda}\|b_\pi\|_{L^1(\mu)}\Big).
\]
\end{lemma}

\begin{proof}
{\sc Step 1.} (Positivity). Since $\lambda$ is the smallest eigenvalue of $\Sigma_T$, the direction of the tilt satisfies
\begin{equation}\label{eq:vbound}
|v(x)|=|\Sigma_T^{-1}b_\pi(x)|\stackrel{\eqref{eq:MT}}{\leq}\lambda^{-1}|b_\pi(x)|\leq\lambda^{-1}\|b_\pi\|_{L^\infty(\mu)}\qquad\text{for }\mu\text{-a.e. }x.
\end{equation}
In particular $v$ is bounded. Since $\pi_x\geq0$, $0\leq s\leq1$ and $\tau_x\leq|\tau_x|$, we obtain for $\mu$-a.e.\ $x$
\[
\tilde\pi_x=(1-s)\,\pi_x+s\,\bar\pi_x-(1-s)\,\tau_x\geq s\,\bar\pi_x-|\tau_x|\stackrel{\eqref{ass:A}}{\geq}s\kappa\,\nu-|\tau_x|\stackrel{\eqref{eq:taubound},\eqref{eq:vbound}}{\geq}s\kappa\,\nu-\frac{4T}{\lambda}\|b_\pi\|_{L^\infty(\mu)}\,\nu\stackrel{\eqref{eq:s}}{\ge}0.
\]

\medskip

{\sc Step 2.} (Mass, mean and marginals). By Lemma~\ref{lem:MT}\ref{it:MT-id}, the tilt has zero mass and first moment $\Sigma_Tv(x)=b_\pi(x)$. Hence $\tilde\pi_x(\mathbb{R}^d)=(1-s)+s-(1-s)\,\tau_x(\mathbb{R}^d)=1$ and
\begin{align*}
\int y\,\tilde\pi_x(dy)&=(1-s)\int y\,\pi_x(dy)+s\int y\,\bar\pi_x(dy)-(1-s)\int y\,\tau_x(dy)\\
&\stackrel{\eqref{eq:defect},\eqref{ass:A}}{=}(1-s)\big(x+b_\pi(x)\big)+s\,x-(1-s)\,\Sigma_Tv(x)=x,
\end{align*}
where we used that $\bar\pi$ is a martingale coupling. The first marginal of $\tilde\pi$ is $\mu$ by construction. For the second marginal, since $\pi$ and $\bar\pi$ have second marginal $\nu$,
\[
\int\tilde\pi_x\,\mu(dx)=(1-s)\,\nu+s\,\nu-(1-s)\Big\langle \Sigma_T^{-1}\int b_\pi\,d\mu,\varphi_T\Big\rangle\nu\stackrel{\eqref{eq:zeromean}}{=}\nu.
\]
Hence $\tilde\pi\in\Pi_M(\mu,\nu)$.

\medskip

{\sc Step 3.} (Weighted distance). By construction $\tilde\pi-\pi=s(\bar\pi-\pi)-(1-s)\,\mu\otimes\tau$. Since $\pi$ and $\bar\pi$ are couplings of $\mu$ and $\nu$, the first term contributes at most
\[
s\int(1+|x|+|y|)\,(\bar\pi+\pi)(dx,dy)\stackrel{\eqref{eq:rAomega}}{\leq}2sA.
\]
For the tilt, using $|x|\leq r$ for $\mu$-a.e.\ $x$, we obtain
\begin{align*}
\int(1+|x|+|y|)\,|\mu\otimes\tau|(dx,dy)&\stackrel{\eqref{eq:taubound}}{\leq}4T\int\!\!\int(1+|x|+|y|)\,|v(x)|\,\nu(dy)\,\mu(dx)\\
&\stackrel{\eqref{eq:vbound}}{\leq}
\frac{4T}{\lambda}
\int|b_\pi(x)|
\Big(1+r+\int|y|\,\nu(dy)\Big)\,\mu(dx)\\
& \stackrel{\eqref{eq:rAomega}}{=}\frac{4TA}{\lambda}\,\|b_\pi\|_{L^1(\mu)}.
\end{align*}
\end{proof}

\subsection{Repair of a martingale bias that is small in \texorpdfstring{$L^1(\mu)$}{L1}}
By Lemma \ref{lem:shadow} the shadow of a martingale coupling has a martingale bias, that is  only small in $L^1(\mu)$, so a few points may still carry a large bias. Since Lemma~\ref{lem:smalldefect} assumes a bias that is uniformly small, we first need to get rid of these points. At every point $x$ with large martingale bias, we therefore replace $\pi_x$ with $\bar\pi_x$, as $\bar\pi$ has no bias at all. This operation has little 
{
\color{blue}
transport
}
cost, since the set of such points has small $\mu$-measure by Markov's inequality. It does however have one side effect: the new coupling no longer has second marginal $\nu$ in general.
Indeed, the mass that used to be sent along $\pi_x$ is now sent along $\bar\pi_x$ for these $x$. The difference is a signed measure $\zeta$ with zero total mass, since the same amount of mass is removed and added. This is where assumption \ref{ass:A} comes into play: we repair this mismatch by mixing in a small amount of $\bar\pi$ and subtracting $\zeta$. This restores the second marginal $\nu$ without creating negative mass, thanks to the lower bound $\bar\pi_x\geq\kappa\nu$. Now that the martingale bias is uniformly small, Lemma~\ref{lem:smalldefect} exactly cancels the remaining martingale bias. 

\begin{proposition}\label{prop:repair}
Let $\omega$ be as in~\eqref{eq:rAomega}, and let
\[
C_1:=4(1+r)+\frac{3A}{\kappa^2}+\frac{16TA}{\kappa\lambda}+\frac{8TA}{\lambda},\qquad \beta_0:=\min\Big\{1,\ \kappa^4,\ \frac{\kappa^2\lambda^2}{64T^2}\Big\}.
\]
If $\pi\in\Pi(\mu,\nu)$ satisfies $\|b_\pi\|_{L^1(\mu)}\leq\beta_0$, there is $\tilde\pi\in\Pi_M(\mu,\nu)$ with
\[
\int(1+|x|+|y|)\,|\tilde\pi-\pi|(dx,dy)\leq C_1\|b_\pi\|_{L^1(\mu)}^{\frac 12 }+4\,\omega\big(\|b_\pi\|_{L^1(\mu)}^{\frac 12}\big).
\]
In particular
\[
\int c\,d\tilde\pi\leq\int c\,d\pi+C\Big(C_1\|b_\pi\|_{L^1(\mu)}^{\frac12}+4\,\omega\big(\|b_\pi\|_{L^1(\mu)}^{\frac12}\big)\Big)
\]
whenever $|c(x,y)|\leq C(1+|x|+|y|)$.
\end{proposition}

\begin{proof}
If $\|b_\pi\|_{L^1(\mu)}=0$, then $\pi$ is already a martingale coupling and we take $\tilde\pi:=\pi$. Therefore, we assume w.l.o.g that $\|b_\pi\|_{L^1(\mu)}>0$. Recall that by assumption $\|b_\pi\|_{L^1(\mu)}\leq\beta_0\leq1$, so that
\begin{equation}\label{eq:normsqrt}
\|b_\pi\|_{L^1(\mu)}\leq\|b_\pi\|_{L^1(\mu)}^{\frac12}.
\end{equation}
Moreover, let us fix a martingale coupling $\bar\pi\in\Pi_M(\mu,\nu)$ satisfying \eqref{ass:A}, so that $\kappa\,\nu\leq\bar\pi_x\leq\kappa^{-1}\nu$ for $\mu$-a.e.\ $x$.

\medskip

{\sc Step 1.} (Points with large bias). Let
\[
E:=\big\{x:|b_\pi(x)|>\|b_\pi\|_{L^1(\mu)}^{\frac12}\big\}.
\]
By Markov's inequality,
\begin{equation}\label{eq:aMarkov}
\mu(E)\leq\frac{1}{\|b_\pi\|_{L^1(\mu)}^{\frac12}}\int|b_\pi|\,d\mu=\|b_\pi\|_{L^1(\mu)}^{\frac12}.
\end{equation}
The quantity $\int_E\pi_x\,\mu(dx)$ is the second marginal of the restriction of $\pi$ to $E\times\mathbb{R}^d$. Therefore it is dominated by $\nu$ and has mass $\mu(E)$. Since $|x|\leq r$ for $\mu$-a.e.\ $x$, we obtain
\begin{equation}\label{eq:Ebound}
\begin{split}
\int_E\int (1+|x|+|y|)\,\pi_x(dy)\,\mu(dx)&=\int_E(1+|x|)\,\mu(dx)+\int|y|\int_E\pi_x(dy)\,\mu(dx)\\
&\stackrel{\eqref{eq:omega}}{\leq}\mu(E)(1+r)+\omega\big(\mu(E)\big).
\end{split}
\end{equation}
The same argument applies to $\bar\pi$, so~\eqref{eq:Ebound} also holds with $\bar\pi$ in place of $\pi$ (where we keep $\pi$ in the definition of $E$).

\medskip

{\sc Step 2.} (Removing large bias). On $E$ we replace the kernel of $\pi$ by the kernel of $\bar\pi$, that is we consider the measure
\[
\pi'(dx,dy):=\mu(dx)\big(\mathbf{1}_{E^c}(x)\,\pi_x(dy)+\mathbf{1}_E(x)\,\bar\pi_x(dy)\big).
\]
By construction its first marginal is $\mu$. Since $\bar\pi$ is a martingale coupling, the bias of $\pi'$ vanishes on $E$, and
\begin{equation}\label{eq:defectprime}
b_{\pi'}(x) = b_\pi(x)\mathbf{1}_{E^c}(x),\qquad\text{so that}\qquad|b_{\pi'}(x)|\leq\|b_\pi\|_{L^1(\mu)}^{\frac 12}
\end{equation}
by definition of $E$. However, the second marginal of $\pi'$ is no longer $\nu$. Indeed,
\[
\int\pi'_x\,\mu(dx)=\int\pi_x\,\mu(dx)+\int_E\big(\bar\pi_x-\pi_x\big)\,\mu(dx)=\nu+\zeta,\qquad\zeta:=\int_E\big(\bar\pi_x-\pi_x\big)\,\mu(dx).
\]
The signed measure $\zeta$ has zero mass, since both terms have mass $\mu(E)$. Moreover
\begin{equation}\label{eq:zetabounds}
|\zeta|\leq\int_E(\bar\pi_x+\pi_x)\,\mu(dx),\qquad\zeta^+\leq\int_E\bar\pi_x\,\mu(dx)\stackrel{\eqref{ass:A}}{\leq}\kappa^{-1}\mu(E)\,\nu.
\end{equation}
Finally, $\pi'-\pi=\mu(dx)\,\mathbf{1}_E(x)\,(\bar\pi_x-\pi_x)(dy)$, so that by \eqref{eq:Ebound}
\begin{equation}\label{eq:step2}
\begin{split}
\int(1+|x|+|y|)\,|\pi'-\pi|(dx,dy)&\leq\int_E\int (1+|x|+|y|)\,(\bar\pi_x+\pi_x)(dy)\,\mu(dx)\\
&\stackrel{\eqref{eq:Ebound}}{\leq}2\mu(E)(1+r)+2\omega\big(\mu(E)\big).
\end{split}
\end{equation}

\medskip

{\sc Step 3.} (Restoring the second marginal). We absorb $\zeta$ by mixing in a small amount of $\bar\pi$ and subtracting the product $\mu\otimes\zeta$. Let $s:=\mu(E)/\kappa^2$. By~\eqref{eq:aMarkov} and $\|b_\pi\|_{L^1(\mu)}\leq\beta_0\leq\kappa^4$ we have $s\leq1$. Define
\[
\pi_0(dx,dy):=\mu(dx)\,\pi_{0,x}(dy),\qquad\pi_{0,x}:=(1-s)\,\pi'_x+s\,\bar\pi_x-(1-s)\,\zeta.
\]
The kernel is nonnegative, since $\pi'_x\geq0$ and $s\leq1$ give
\[
\pi_{0,x}\geq s\,\bar\pi_x-\zeta^+\stackrel{\eqref{ass:A},\eqref{eq:zetabounds}}{\geq}s\kappa\,\nu-\kappa^{-1}\mu(E)\,\nu=0.
\]
It has mass $(1-s)+s-(1-s)\,\zeta(\mathbb{R}^d)=1$. The first marginal of $\pi_0$ is $\mu$ and its second marginal is
\[
\int\pi_{0,x}\,\mu(dx)=(1-s)(\nu+\zeta)+s\,\nu-(1-s)\,\zeta=\nu,
\]
so that $\pi_0\in\Pi(\mu,\nu)$. To compute its bias, note that
\[
\int y\,\zeta(dy)=\int_E\Big(\int y\,\bar\pi_x(dy)-\int y\,\pi_x(dy)\Big)\,\mu(dx)=-\int_E b_\pi\,d\mu,
\]
where we used that $\bar \pi$ is a martingale coupling. Hence, using again that $\bar\pi_x$ has mean $x$,
\[
b_{\pi_0}(x)=(1-s)\,b_{\pi'}(x)-(1-s)\int y\,\zeta(dy)\stackrel{\eqref{eq:defectprime}}{=}(1-s)\Big(b_\pi(x)\mathbf{1}_{E^c}(x)+\int_E b_\pi\,d\mu\Big).
\]
The first term is bounded by $\|b_\pi\|_{L^1(\mu)}^{1/2}$ and the second by $\|b_\pi\|_{L^1(\mu)}$, so that
\begin{equation}\label{eq:defectzero}
\|b_{\pi_0}\|_{L^\infty(\mu)}\stackrel{\eqref{eq:normsqrt}}{\leq}2\|b_\pi\|_{L^1(\mu)}^{\frac 12},\qquad\|b_{\pi_0}\|_{L^1(\mu)}\leq2\|b_\pi\|_{L^1(\mu)}.
\end{equation}
For the transport cost, note that by definition $\pi_0-\pi'=s(\bar\pi-\pi')-(1-s)\,\mu\otimes\zeta$. Since $\pi'\leq\pi+\bar\pi$ and every coupling of $\mu$ and $\nu$ integrates $1+|x|+|y|$ to at most $A$,
\begin{align*}
s\int(1+|x|+|y|)\,(\bar\pi+\pi')(dx,dy)&\leq2s\int(1+|x|+|y|)\,\bar\pi(dx,dy)\\
&\quad+s\int(1+|x|+|y|)\,\pi(dx,dy)\leq3As.
\end{align*}
For the product term, $|\zeta|$ is dominated by the sum of two measures which are dominated by $\nu$ and have mass $\mu(E)$, see~\eqref{eq:zetabounds}. Hence
\begin{align*}
\int\!\!\int (1+|x|+|y|)\,\mu(dx)\,|\zeta|(dy)&=\Big(1+\int|x|\,d\mu\Big)\|\zeta\|_{TV}+\int|y|\,|\zeta|(dy)\\
&\stackrel{\eqref{eq:omega}}{\leq}2\mu(E)(1+r)+2\omega\big(\mu(E)\big),
\end{align*}
and therefore
\begin{equation}\label{eq:step3}
\int(1+|x|+|y|)\,|\pi_0-\pi'|(dx,dy)\leq3As+2\mu(E)(1+r)+2\omega\big(\mu(E)\big).
\end{equation}

\medskip

{\sc Step 4.} (Exact cancellation). The bias of $\pi_0$ is uniformly small. Indeed, by the choice of $\beta_0$,
\[
\|b_{\pi_0}\|_{L^\infty(\mu)}\stackrel{\eqref{eq:defectzero}}{\leq}2\|b_\pi\|_{L^1(\mu)}^{1/2}\leq2\beta_0^{1/2}\leq\frac{\kappa\lambda}{4T}.
\]
Therefore we may apply Lemma~\ref{lem:smalldefect} to $\pi_0\in\Pi(\mu,\nu)$ and obtain $\tilde\pi\in\Pi_M(\mu,\nu)$ with parameter
\[
s'=\frac{4T\|b_{\pi_0}\|_{L^\infty(\mu)}}{\kappa\lambda}\stackrel{\eqref{eq:defectzero}}{\leq}\frac{8T}{\kappa\lambda}\|b_\pi\|_{L^1(\mu)}^{\frac 12},
\]
and
\begin{equation}\label{eq:step4}
\int(1+|x|+|y|)\,|\tilde\pi-\pi_0|(dx,dy)\leq A\Big(2s'+\frac{4T}{\lambda}\|b_{\pi_0}\|_{L^1(\mu)}\Big)\stackrel{\eqref{eq:defectzero}}{\leq}A\Big(\frac{16T}{\kappa\lambda}\|b_\pi\|_{L^1(\mu)}^{\frac 12}+\frac{8T}{\lambda}\|b_\pi\|_{L^1(\mu)}\Big).
\end{equation}

\medskip

{\sc Step 5.} (Conclusion). By the triangle inequality,
\begin{align*}
\int(1+|x|+|y|)\,|\tilde\pi-\pi|(dx,dy)\leq{}&\int(1+|x|+|y|)\,|\tilde\pi-\pi_0|(dx,dy)\\
&+\int(1+|x|+|y|)\,|\pi_0-\pi'|(dx,dy)\\
&+\int(1+|x|+|y|)\,|\pi'-\pi|(dx,dy).
\end{align*}
Combining the latter with~\eqref{eq:step4}, \eqref{eq:step3} and~\eqref{eq:step2}, and using~\eqref{eq:aMarkov}, \eqref{eq:normsqrt} and the monotonicity of $\omega$, we get
\[
\int(1+|x|+|y|)\,|\tilde\pi-\pi|(dx,dy)\leq\Big(4(1+r)+\frac{3A}{\kappa^2}+\frac{16TA}{\kappa\lambda}+\frac{8TA}{\lambda}\Big)\|b_\pi\|_{L^1(\mu)}^{\frac 12}+4\,\omega\big(\|b_\pi\|_{L^1(\mu)}^{\frac 12}\big).
\]
The bracket on the right-hand side is exactly $C_1$, which proves the first claim. For the second claim, if $|c(x,y)|\leq C(1+|x|+|y|)$, then
\[
\Big|\int c\,d\tilde\pi-\int c\,d\pi\Big|\leq\int|c|\,d|\tilde\pi-\pi|\leq C\int(1+|x|+|y|)\,|\tilde\pi-\pi|(dx,dy),
\]
which concludes the proof.
\end{proof}

\subsection{Uniform penalization}
In the previous subsection we showed that a coupling with a small martingale bias can be turned into a martingale coupling at a small transport cost. We now use this observation to compare $V$ and $P_R$. Since martingale couplings have a bias of zero, we always have $P_R\leq V$. It thus remains to bound $V-P_R$ from above. 
For this, given an arbitrary $\pi\in\Pi(\mu,\nu)$, we write
$\beta:=\|b_\pi\|_{L^1(\mu)}$. If $\beta$ is sufficiently small,
Proposition~\ref{prop:repair} produces a martingale coupling with
a change in transport cost bounded by
$C\bigl(C_1\sqrt{\beta}+4\omega(\sqrt{\beta})\bigr)$.
This repair error tends to zero as $\beta\to0$. For biases bounded
away from zero, the penalty $R\beta$ dominates
the repair error when $R$ is large. Thus, the largest repair error
not covered by the penalty tends to zero as $R\to\infty$.
If instead $\beta$ exceeds the threshold for applying
Proposition~\ref{prop:repair}, we compare $\pi$ directly with the
reference martingale coupling $\bar\pi$ satisfying~\eqref{ass:A}.
The linear-growth bound on $c$ then ensures that, for sufficiently
large $R$, the penalized cost of $\pi$ is at least the cost of
$\bar\pi$. Together, these two cases give an upper bound on
$V-P_R$ that tends to zero as $R\to\infty$. The constants and error
bounds depend only on $\kappa,T,\lambda,r,A$, the tail modulus
$\omega$, and the growth constant of $c$. By
Lemma~\ref{lem:uniformconst}, these data can be chosen uniformly
for all sufficiently large $n$, giving the uniform comparison
needed to prove Theorem~\ref{thm:main}.

\begin{corollary}\label{cor:uniform}
Let $|c(x,y)|\leq C(1+|x|+|y|)$. Let $C_1$ and $\beta_0$ be as in Proposition~\ref{prop:repair}, and set
\[
\Psi(\beta):=C\Big(C_1\beta^{\frac12}+4\,\omega\big(\beta^{\frac12}\big)\Big),\qquad R_*:=\frac{2CA}{\beta_0},
\]
and
\begin{equation}\label{eq:eR}
e(R):=\sup_{0\leq\beta\leq\beta_0}\big(\Psi(\beta)-R\beta\big),\qquad R\geq0.
\end{equation}
Then, for every $R\geq R_*$,
\begin{equation}\label{eq:VPR}
0\leq V(\mu,\nu)-P_R(\mu,\nu)\leq e(R),
\end{equation}
and $e(R)\to0$ as $R\to\infty$. If moreover $\omega(a)\leq Ka^{1-1/p}$ for some $p\in(1,\infty)$ and $K\geq0$, then, with $\gamma:=\frac{1-1/p}{2}$ and $C_2:=C(C_1+4K)$,
\begin{equation}\label{eq:eRrate}
e(R)\leq C_2^{\frac1{1-\gamma}}R^{-\frac\gamma{1-\gamma}}.
\end{equation}
The same holds with $\gamma=\frac12$ if $\omega(a)\leq Ka$.
\end{corollary}

\begin{proof}
{\sc Step 1.} (Comparison between $V$ and $P_R$). We first show the lower bound in~\eqref{eq:VPR}. If $\pi\in\Pi_M(\mu,\nu)$, then $b_\pi=0$ for $\mu$-a.e.\ $x$, so that the penalized cost of $\pi$ equals $\int c\,d\pi$. Since $\Pi_M(\mu,\nu)\subset\Pi(\mu,\nu)$, the infimum in~\eqref{eq:PR} is taken over a larger set than the one defining $V$, and $P_R(\mu,\nu)\leq V(\mu,\nu)$.

\medskip

For the upper bound in~\eqref{eq:VPR}, by the definition~\eqref{eq:PR} of $P_R$  it suffices to show that every $\pi\in\Pi(\mu,\nu)$ satisfies
\begin{equation}\label{eq:VPRpi}
V(\mu,\nu)-\int c\,d\pi-R\,\|b_\pi\|_{L^1(\mu)}\leq e(R).
\end{equation}
We distinguish between two cases. If $\|b_\pi\|_{L^1(\mu)}\leq\beta_0$, Proposition~\ref{prop:repair} gives $\tilde\pi\in\Pi_M(\mu,\nu)$ with $\int c\,d\tilde\pi\leq\int c\,d\pi+\Psi(\|b_\pi\|_{L^1(\mu)})$. Since $\tilde\pi$ is admissible for $V$,
\begin{align*}
V(\mu,\nu)-\int c\,d\pi-R\,\|b_\pi\|_{L^1(\mu)}&\leq\int c\,d\tilde\pi-\int c\,d\pi-R\,\|b_\pi\|_{L^1(\mu)}\\
&\leq\Psi\big(\|b_\pi\|_{L^1(\mu)}\big)-R\,\|b_\pi\|_{L^1(\mu)}\stackrel{\eqref{eq:eR}}{\leq}e(R).
\end{align*}
If $\|b_\pi\|_{L^1(\mu)}>\beta_0$, we compare $\pi$ with the martingale coupling $\bar\pi$ satisfying \eqref{ass:A}. Every coupling of $\mu$ and $\nu$ integrates $1+|x|+|y|$ to at most $A$, see~\eqref{eq:rAomega}. Hence $V(\mu,\nu)\leq\int c\,d\bar\pi\leq CA$ and $\int c\,d\pi\geq-CA$, so that for $R\geq R_*$
\[
V(\mu,\nu)-\int c\,d\pi-R\,\|b_\pi\|_{L^1(\mu)}\leq2CA-R\beta_0\leq0\leq\Psi(0)\stackrel{\eqref{eq:eR}}{\leq}e(R),
\]
where $\Psi(0)=0$ follows from $\omega(0)=0$. This proves~\eqref{eq:VPRpi}, and taking the infimum over $\pi$ gives the upper bound in~\eqref{eq:VPR}.

\medskip

{\sc Step 2.} (Decay of $e$). The function $\Psi$ is nondecreasing and $\Psi(\beta)\to0$ as $\beta\to0$, since $\omega(a)\to0$ as $a\to0$ by~\eqref{eq:omegaUI}. Fix $\delta\in(0,\beta_0]$ and split the supremum in~\eqref{eq:eR}. For $\beta\leq\delta$ we have $\Psi(\beta)-R\beta\leq\Psi(\delta)$. For $\beta\in(\delta,\beta_0]$ we have
\[
\Psi(\beta)-R\beta\leq\Psi(\beta_0)-R\delta\leq0\qquad\text{as soon as }R\geq\Psi(\beta_0)/\delta.
\]
Hence $e(R)\leq\Psi(\delta)$ for $R$ large, so that $\limsup_{R\to\infty}e(R)\leq\Psi(\delta)$. Letting $\delta\to0$ gives $e(R)\to0$.

\medskip

{\sc Step 3.} (Rate). Let $0\leq\beta\leq\beta_0\leq1$. Since $\gamma\leq 1/2$, we have $\beta^{1/2}\leq\beta^\gamma$, and the assumption on $\omega$ gives $\omega(\beta^{1/2})\leq K\beta^{(1-1/p)/2}=K\beta^\gamma$. Therefore $\Psi(\beta)\leq C(C_1+4K)\beta^\gamma=C_2\beta^\gamma$ and
\[
e(R)\leq\sup_{\beta\geq0}\big(C_2\beta^\gamma-R\beta\big).
\]
The supremum is attained at $\beta_*:=(\gamma C_2/R)^{1/(1-\gamma)}$, where $R\beta_*=\gamma C_2\beta_*^\gamma$. Its value is
\[
(1-\gamma)\,C_2\,\beta_*^\gamma=(1-\gamma)\gamma^{\frac \gamma{1-\gamma}}C_2^{\frac 1{1-\gamma}}R^{-\frac \gamma{1-\gamma}}\leq C_2^{\frac 1{1-\gamma}}R^{- \frac \gamma{1-\gamma}},
\]
since the prefactor is at most one. The latter is~\eqref{eq:eRrate}. Finally, if $\omega(a)\leq Ka$, the same argument with $\gamma=1/2$ applies, since then $\omega(\beta^{1/2})\leq K\beta^{1/2}$.
\end{proof}
\section{Proof of Theorem~\ref{thm:main}}\label{sec:proof}

In this section we conclude the proof of Theorem~\ref{thm:main}. The constants appearing in Proposition~\ref{prop:repair} and Corollary~\ref{cor:uniform} depend on the pair $(\mu,\nu)$, and we need them to be uniform along the sequence $(\mu_n,\nu_n)$. This is the content of Lemma~\ref{lem:uniformconst}, and it is where the non-degeneracy of $\nu$ and the Wasserstein convergence of $(\mu_n,\nu_n)$ are used. We then combine the uniform comparison between $V$ and $P_R$ of Corollary~\ref{cor:uniform} with the stability of $P_R$ of Proposition~\ref{prop:PRstab}. Finally, we choose $R$ depending on $n$ to obtain a convergence rate.

\begin{lemma}\label{lem:uniformconst}
Under the hypotheses of Theorem~\ref{thm:main} there are $T>0$, $\lambda>0$, $r,A<\infty$, a nondecreasing function $\omega$ with $\omega(a)\to0$ as $a\to0$, and $n_0\in\mathbb{N}$ such that
\begin{equation}\label{eq:uniformconst}
\lambda_T(\nu_n)\geq\lambda,\qquad
\operatorname*{ess\,sup}_{\mu_n}|x|\leq r,\qquad
1+r+\int|y|\,\nu_n(dy)\leq A,\qquad
\omega_{\nu_n}\leq\omega
\end{equation}
for every $n\geq n_0$, and the same holds with $(\mu,\nu)$ in place of $(\mu_n,\nu_n)$. In particular, Corollary~\ref{cor:uniform} holds for all these pairs with the same $R_*$ and the same function $e$. Moreover, under the moment assumption of Theorem~\ref{thm:main}\ref{it:main-rate} one may choose $\omega(a)=Ka^{1-1/p}$ for some constant $K$. If the supports of all $\nu_n$ lie in a fixed ball $B(0,K)$, one may choose $\omega(a)=Ka$.
\end{lemma}

\begin{proof}
We start by showing that the measure $\nu_n$ satisfies Assumption \ref{ass:2} for large enough $n$. By Lemma~\ref{lem:Alimit}, the limit pair $(\mu,\nu)$ satisfies Assumption \ref{ass:1} with the same $\kappa$. Since $\nu$ is not carried by a hyperplane, Lemma~\ref{lem:MT}\ref{it:MT-pos} gives $T>0$ with $\lambda_T(\nu)>0$, and we set $\lambda:=\lambda_T(\nu)/2$. By Lemma~\ref{lem:MT}\ref{it:MT-cont}, $\lambda_T(\nu_n)\to\lambda_T(\nu)=2\lambda$, so that there is $n_0$ with $\lambda_T(\nu_n)\geq\lambda$ for every $n\geq n_0$. This is the only place where the non-degeneracy of $\nu$ is used.

\medskip

Next, we bound the supports of $\mu_n$ and the first moments of $\nu_n$. Since $\nu_n\to\nu$ in $W_1$, the quantities $m(\nu_n)$, $\int|y|\,\nu_n(dy)$ and $\widehat M_1(\nu_n)$ converge, and are therefore bounded in $n$. Since every pair $(\mu_n,\nu_n)$ satisfies Assumption \ref{ass:1}, Lemma~\ref{lem:support} gives
\[
\operatorname*{ess\,sup}_{\mu_n}|x|\leq|m(\nu_n)|+\kappa^{-1}\widehat M_1(\nu_n)\leq\sup_k\Big(|m(\nu_k)|+\kappa^{-1}\widehat M_1(\nu_k)\Big)<\infty,
\]
and the same bound holds for $\mu$, since $(\mu,\nu)$ satisfies Assumption \ref{ass:1} as well. We let $r$ be the maximum of these bounds, and define $A:=1+r+\sup_n\int|y|\,\nu_n(dy)\vee\int|y|\,\nu(dy)$.

\medskip

For the tails, we set $\omega:=\sup_n\omega_{\nu_n}\vee\omega_\nu$. As a supremum of nondecreasing functions, $\omega$ is nondecreasing, and $\omega(a)\to0$ as $a\to0$ by~\eqref{eq:omegaUI}.

\medskip

Finally, we consider the two rate assumptions. Under the moment assumption, let $K:=\sup_n M_p(\nu_n)^{1/p}$. By lower semicontinuity of $\nu\mapsto\int|y|^p\,d\nu$ under weak convergence, $M_p(\nu)\leq\liminf_n M_p(\nu_n)\leq K^p$, so that~\eqref{eq:omegaHolder} gives
\[
\omega_{\nu_n}(a)\leq Ka^{1-\frac1 p}\quad\text{for all }n,\qquad\omega_\nu(a)\leq Ka^{1-\frac 1 p},
\]
and we may take $\omega(a)=Ka^{1-1/p}$. If instead the supports of all $\nu_n$ lie in $B(0,K)$, so does the support of the weak limit $\nu$, and $\omega_{\nu_n}(a),\omega_\nu(a)\leq Ka$ by~\eqref{eq:omega}, so that we may take $\omega(a)=Ka$.
\end{proof}

\begin{proof}[Proof of Theorem~\ref{thm:main}]
Throughout the proof we fix $T$, $\lambda$, $r$, $A$, $\omega$ and $n_0$ as in Lemma~\ref{lem:uniformconst}, and we let $C_1$, $\beta_0$, $R_*$ and $e$ be the corresponding quantities of Proposition~\ref{prop:repair} and Corollary~\ref{cor:uniform}. By Lemma~\ref{lem:uniformconst}, none of them depends on $n$, and Corollary~\ref{cor:uniform} holds with the same $R_*$ and the same $e$ for every pair $(\mu_n,\nu_n)$ with $n\geq n_0$ and for the limit pair $(\mu,\nu)$.

\medskip

{\sc Step 1.} (Proof of~\ref{it:main-cont}). The lower bound $\liminf_n V(\mu_n,\nu_n)\geq V(\mu,\nu)$ holds by Corollary~\ref{cor:lsc}, without any assumption on the pairs. For the upper bound, we fix $R\geq R_*$. We first substitute $V(\mu_n,\nu_n)$ by $P_R(\mu_n,\nu_n)$, which costs at most $e(R)$ uniformly in $n$. Then we pass to the limit in $P_R$, and finally we go back from $P_R$ to $V$ at the limit pair. Summing up this yields
\[
\limsup_{n\to\infty}V(\mu_n,\nu_n)\stackrel{\eqref{eq:VPR}}{\leq}\limsup_{n\to\infty}P_R(\mu_n,\nu_n)+e(R)\leq P_R(\mu,\nu)+e(R)\stackrel{\eqref{eq:VPR}}{\leq}V(\mu,\nu)+e(R),
\]
where the second inequality follows from Proposition~\ref{prop:PRstab}\ref{it:PR-cont}. Since the left-hand side does not depend on $R$, we may let $R\to\infty$. By Corollary~\ref{cor:uniform}, $e(R)\to0$, so that $\limsup_n V(\mu_n,\nu_n)\leq V(\mu,\nu)$. Together with the lower bound, this proves~\ref{it:main-cont}.

\medskip

{\sc Step 2.} (A two-sided estimate). We now work under the assumptions of~\ref{it:main-rate}. Since $c$ is $L$-Lipschitz, we have $|c(x,y)|\leq|c(0,0)|+L(|x|+|y|)\leq C(1+|x|+|y|)$, so that we are in the setting of Corollary~\ref{cor:uniform}. Moreover, by Lemma~\ref{lem:uniformconst} we may take $\omega(a)=Ka^{1-1/p}$, so that~\eqref{eq:eRrate} holds. Fix $R\geq R_*$ and $n\geq n_0$. We argue as in {\sc Step~1}, using the quantitative stability of $P_R$, and we obtain
\[
\begin{split}
V(\mu_n,\nu_n)& \stackrel{\eqref{eq:VPR}}{\leq}P_R(\mu_n,\nu_n)+e(R) \\
& \leq P_R(\mu,\nu)+(L+R)\Delta_n+e(R) \\
& \stackrel{\eqref{eq:VPR}}{\leq}V(\mu,\nu)+(L+R)\Delta_n+e(R),
\end{split}
\]
where the second inequality follows from Proposition~\ref{prop:PRstab}\ref{it:PR-lip}. Exchanging the roles of the two pairs
 and combining the two estimates with~\eqref{eq:eRrate}, we get
\begin{equation}\label{eq:twosided}
\big|V(\mu_n,\nu_n)-V(\mu,\nu)\big|\leq L\Delta_n+R\Delta_n+C_2^{\frac1{1-\gamma}}R^{-\frac\gamma{1-\gamma}}\qquad\text{for all }R\geq R_*.
\end{equation}

\medskip

{\sc Step 3.} (Choice of $R$). In~\eqref{eq:twosided}, the term $R\Delta_n$ is the price of a large penalty and grows with $R$, while the last term is the price of a small penalty and decreases with $R$. If $\Delta_n=0$, letting $R\to\infty$ in~\eqref{eq:twosided} gives $V(\mu_n,\nu_n)=V(\mu,\nu)$. Otherwise, we balance the two terms by choosing
\[
R:=\Delta_n^{-(1-\gamma)},\qquad\text{so that}\qquad R\Delta_n=\Delta_n^\gamma\quad\text{and}\quad R^{-\frac\gamma{1-\gamma}}=\Delta_n^\gamma.
\]
This choice is admissible as soon as $\Delta_n^{-(1-\gamma)}\geq R_*$, which holds for $n$ large since $\Delta_n\to0$ and $\gamma<1$. For such $n$, \eqref{eq:twosided} becomes
\[
\big|V(\mu_n,\nu_n)-V(\mu,\nu)\big|\leq L\Delta_n+\big(1+C_2^{\frac1{1-\gamma}}\big)\Delta_n^\gamma.
\]
The constant $1+C_2^{1/(1-\gamma)}$ does not depend on $n$, which proves~\ref{it:main-rate}. Finally, if the supports of all $\nu_n$ lie in a fixed ball, Lemma~\ref{lem:uniformconst} allows us to take $\omega(a)=Ka$, so that~\eqref{eq:eRrate} holds with $\gamma=\frac12$, and the same argument applies.
\end{proof}

\begin{remark}[Approximation of martingale couplings]\label{rem:approx}
The proof of Theorem~\ref{thm:main} also shows that every martingale coupling of the limit can be approximated by martingale couplings along the sequence. This is the converse of Lemma~\ref{lem:compact}\ref{it:compact-mart}, which fails in general, as shown in   \cite{BJ22}. Indeed, let $\pi\in\Pi_M(\mu,\nu)$ and let $\tilde\pi_n$ be its shadow on $(\mu_n,\nu_n)$. Since $b_\pi=0$, Lemma~\ref{lem:shadow}\ref{it:shadow-defect} gives $\|b_{\tilde\pi_n}\|_{L^1(\mu_n)}\leq\Delta_n$, which is at most $\beta_0$ for $n$ large. Applying Proposition~\ref{prop:repair} to $\tilde\pi_n$, with the constants of Lemma~\ref{lem:uniformconst}, we obtain $\pi_n\in\Pi_M(\mu_n,\nu_n)$ with
\[
\int(1+|x|+|y|)\,|\pi_n-\tilde\pi_n|(dx,dy)\leq C_1\Delta_n^{\frac12}+4\,\omega\big(\Delta_n^{\frac12}\big).
\]
Since the left-hand side bounds $W_1(\pi_n,\tilde\pi_n)$, see \cite[Theorem~6.15]{Vil09}, and $W_1(\tilde\pi_n,\pi)\leq\Delta_n$ by Lemma~\ref{lem:shadow}\ref{it:shadow-dist}, the triangle inequality gives
\[
W_1(\pi_n,\pi)\leq\Delta_n+C_1\Delta_n^{\frac12}+4\,\omega\big(\Delta_n^{\frac12}\big)\to0.
\]
Under the moment assumption of Theorem~\ref{thm:main}\ref{it:main-rate}, Lemma~\ref{lem:uniformconst} allows $\omega(a)=Ka^{1-1/p}$, and the right-hand side is of order $\Delta_n^\gamma$.
\end{remark}

\section{\texorpdfstring{Assumption~\eqref{ass:A} alone is not sufficient}{Assumption (A kappa) alone is not sufficient}}\label{sec:example}

In the counterexample of \cite{BJ22}, every point can send its mass only to a few prescribed points of the target. In particular, the marginals there cannot satisfy~\eqref{ass:A}, which requires the kernels of one martingale coupling to charge the whole target. It is therefore natural to ask whether~\eqref{ass:A} alone, with a constant $\kappa$ independent of $n$, is enough for the continuity of $V$. The following example shows that this is not the case. In particular, the non-degeneracy of $\nu$ in Theorem~\ref{thm:main} cannot be dropped.

\begin{example}\label{ex:parabola}
Let $U_r$ denote the uniform law on $[-r,r]$, and let $\varepsilon_n:=1/n$ and $t_n:=1/n^2$. Let $X_1,W\sim U_{1/4}$ and $Y_1,Z\sim U_1$ be independent, and set
\[
\mu_n:=\operatorname{law}\Big(X_1,\frac{\varepsilon_n}{3}+t_nW\Big),\qquad \nu_n:=\operatorname{law}\big(Y_1,\varepsilon_nY_1^2+t_nZ\big),\qquad c(x,y):=(1+x_1)(y_1-x_1)^2.
\]
Then $\mu_n\to\mu:=U_{1/4}\otimes\delta_0$ and $\nu_n\to\nu:=U_1\otimes\delta_0$ in $W_1$, and every pair $(\mu_n,\nu_n)$ satisfies~\eqref{ass:A} with $\kappa=1/16$. However,
\[
\lim_{n\to\infty}V(\mu_n,\nu_n)\geq V(\mu,\nu)+\frac{1}{540}.
\]
\end{example}

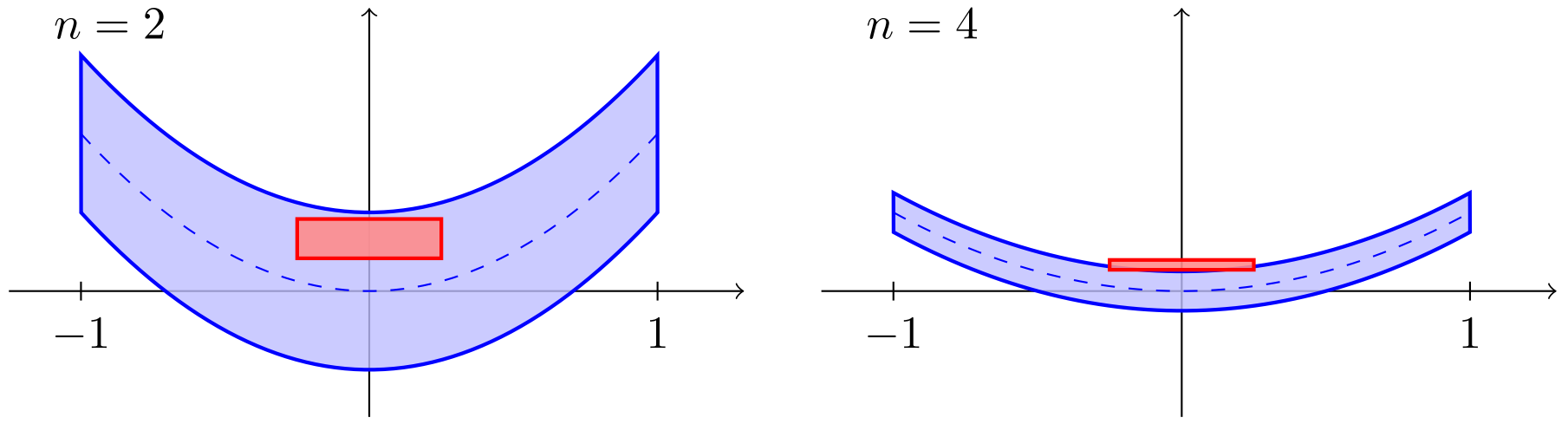
\begin{figure}[ht]
\centering
\begin{tikzpicture}[x=2.2cm,y=2.4cm]
  \begin{scope}
    \draw[->] (-1.25,0) -- (1.3,0);
    \draw[->] (0,-0.4) -- (0,0.9);
    \fill[blue!25, fill opacity=0.8] plot[domain=-1:1, samples=60] (\x,{0.5*\x*\x+0.25}) -- plot[domain=1:-1, samples=60] (\x,{0.5*\x*\x-0.25}) -- cycle;
    \draw[blue, thick] plot[domain=-1:1, samples=60] (\x,{0.5*\x*\x+0.25}) -- plot[domain=1:-1, samples=60] (\x,{0.5*\x*\x-0.25}) -- cycle;
    \draw[blue, dashed] plot[domain=-1:1, samples=60] (\x,{0.5*\x*\x});
    \fill[red!45, fill opacity=0.9] (-0.25,{1/6-1/16}) rectangle (0.25,{1/6+1/16});
    \draw[red, thick] (-0.25,{1/6-1/16}) rectangle (0.25,{1/6+1/16});
    \draw (-1,0.03) -- (-1,-0.03) node[below, font=\small] {$-1$};
    \draw (1,0.03) -- (1,-0.03) node[below, font=\small] {$1$};
    \node[font=\small] at (-0.9,0.85) {$n=2$};
  \end{scope}
  \begin{scope}[xshift=6.2cm]
    \draw[->] (-1.25,0) -- (1.3,0);
    \draw[->] (0,-0.4) -- (0,0.9);
    \fill[blue!25, fill opacity=0.8] plot[domain=-1:1, samples=60] (\x,{0.25*\x*\x+1/16}) -- plot[domain=1:-1, samples=60] (\x,{0.25*\x*\x-1/16}) -- cycle;
    \draw[blue, thick] plot[domain=-1:1, samples=60] (\x,{0.25*\x*\x+1/16}) -- plot[domain=1:-1, samples=60] (\x,{0.25*\x*\x-1/16}) -- cycle;
    \draw[blue, dashed] plot[domain=-1:1, samples=60] (\x,{0.25*\x*\x});
    \fill[red!45, fill opacity=0.9] (-0.25,{1/12-1/64}) rectangle (0.25,{1/12+1/64});
    \draw[red, thick] (-0.25,{1/12-1/64}) rectangle (0.25,{1/12+1/64});
    \draw (-1,0.03) -- (-1,-0.03) node[below, font=\small] {$-1$};
    \draw (1,0.03) -- (1,-0.03) node[below, font=\small] {$1$};
    \node[font=\small] at (-0.9,0.85) {$n=4$};
  \end{scope}
\end{tikzpicture}
\caption{The supports of $\mu_n$ (red) and $\nu_n$ (blue) in Example~\ref{ex:parabola} for $n=2$ and $n=4$, with the vertical axis stretched. As $n\to\infty$, the band around the parabola $\{y_2=\varepsilon_ny_1^2\}$ (dashed) and the support of $\mu_n$ collapse onto the horizontal axis.}
\label{fig:parabola}
\end{figure}

\medskip

In the example, $\nu_n$ lives on a thin band around the parabola $\{y_2=\varepsilon_ny_1^2\}$, which collapses onto a segment of the horizontal axis as $n\to\infty$, see Figure~\ref{fig:parabola}. Every pair $(\mu_n,\nu_n)$ satisfies Assumption~\ref{ass:1} with the same $\kappa=1/16\in (0,1)$ as we will show below, so that every point can send its mass over the whole target. However, since the band is thin, balancing the barycenter in the vertical direction forces every kernel to have almost the same second moment in the horizontal direction. This hidden constraint fixes the cost of all martingale couplings of $(\mu_n,\nu_n)$, up to a small error. In the limit the band becomes a segment, the constraint disappears, and the cost can be lowered.

\medskip

We now verify the claims of Example~\ref{ex:parabola}.
The convergence in $W_1$ holds since $\varepsilon_n,t_n\to0$
and all the random variables involved are bounded.
Moreover, all the supports are contained in
$K:=[-1,1]\times[-2,2]$.
Let $P_K$ denote the orthogonal projection onto $K$, and define
\[
\widehat c(x,y):=c(P_Kx,P_Ky).
\]
Since $c$ is Lipschitz on $K\times K$ and $P_K$ is
$1$-Lipschitz, $\widehat c$ is bounded and globally Lipschitz.
It agrees with $c$ on $K\times K$, so replacing $c$ by
$\widehat c$ leaves all the values $V(\mu_n,\nu_n)$ and
$V(\mu,\nu)$ unchanged.
Hence all the assumptions of Theorem~\ref{thm:main} are satisfied,
except for the non-degeneracy of $\nu$, which is carried by
the line $\{y_2=0\}$.

\medskip
We first construct a martingale coupling $\bar\pi$ satisfying~\eqref{ass:A}. Let $\Phi(x_1,w):=(x_1,\varepsilon_n/3+t_nw)$ and $\Psi(y_1,z):=(y_1,\varepsilon_ny_1^2+t_nz)$. Both maps are injective, and $\mu_n=\Phi_\#(U_{1/4}\otimes U_{1/4})$ and $\nu_n=\Psi_\#(U_1\otimes U_1)$. Let $\bar\pi_n$ be the image under $\Phi\times\Psi$ of
\[
(1+3x_1y_1)(1+3wz)\,U_{1/4}(dx_1)\,U_{1/4}(dw)\,U_1(dy_1)\,U_1(dz).
\]
Integrating out $(x_1,w)$ removes the linear terms, since $U_{1/4}$ is centered, so that the second marginal of $\bar\pi_n$ is $\nu_n$. In the same way, integrating out $(y_1,z)$ shows that its first marginal is $\mu_n$. To check the martingale property, fix $x=\Phi(x_1,w)$. Under the kernel of $\bar\pi_n$ at $x$, the variables $y_1$ and $z$ are independent, and since $\int y^2\,U_1(dy)=1/3$,
\[
\int y_1(1+3x_1y_1)\,U_1(dy_1)=x_1,\qquad\int y_1^2(1+3x_1y_1)\,U_1(dy_1)=\frac{1}{3},\qquad\int z(1+3wz)\,U_1(dz)=w.
\]
Hence the kernel has mean $(x_1,\varepsilon_n/3+t_nw)=x$, and $\bar\pi_n$ is a martingale coupling. Finally, since $\Phi\times\Psi$ is injective, the density of $\bar\pi_n$ with respect to $\mu_n\otimes\nu_n$ is $(1+3x_1y_1)(1+3wz)$, which takes values in $[1/16,49/16]$. Therefore~\eqref{ass:A} holds with $\kappa=1/16$.

\medskip

Next, we show that all martingale couplings of $(\mu_n,\nu_n)$ have almost the same cost. Let $\pi\in\Pi_M(\mu_n,\nu_n)$, let $(X,Y)$ have law $\pi$, and write $X=\Phi(x_1,w)$. On $\operatorname{supp}\nu_n$ the variable $Z:=(Y_2-\varepsilon_nY_1^2)/t_n$ takes values in $[-1,1]$, and the martingale property in the second coordinate reads
\[
\varepsilon_n\mathbb{E}[Y_1^2\mid X]+t_n\mathbb{E}[Z\mid X]=\frac{\varepsilon_n}{3}+t_nw.
\]
Since $|w|\leq1/4$ and $|Z|\leq1$, this gives
\[
\Big|\mathbb{E}[Y_1^2\mid X]-\frac{1}{3}\Big|=\frac{t_n}{\varepsilon_n}\big|w-\mathbb{E}[Z\mid X]\big|\leq\frac{5}{4n}.
\]
This is the hidden constraint: the conditional second moment of $Y_1$ is almost fixed, whatever the martingale coupling. On the other hand, the martingale property in the first coordinate gives $\mathbb{E}[Y_1\mid X]=x_1$, so that
\[
\int c\,d\pi=\int(1+x_1)\big(\mathbb{E}[Y_1^2\mid X]-x_1^2\big)\,\mu_n(dx).
\]
Since the first coordinate of $\mu_n$ has law $U_{1/4}$ and $|1+x_1|\leq5/4$, every martingale coupling of $(\mu_n,\nu_n)$ has cost within $25/(16n)$ of
\[
V_Q:=\int(1+x)\Big(\frac{1}{3}-x^2\Big)\,U_{1/4}(dx).
\]
In particular, $V(\mu_n,\nu_n)\to V_Q$.

\medskip

Finally, we show that the limit problem does better. Since $\mu$ and $\nu$ are carried by the line $\{y_2=0\}$, the martingale couplings of $(\mu,\nu)$ are the one-dimensional martingale couplings of $(U_{1/4},U_1)$, embedded in the first coordinate. For $x\in[-1/4,1/4]$ we consider the kernel
\[
\Big(1+3xy-x\Big(y^2-\frac{1}{3}\Big)\Big)\,U_1(dy).
\]
Its density is at least $1-3/4-1/6>0$. Compared with the kernel $(1+3xy)\,U_1(dy)$ used above, the added term has zero mass and zero mean, since $y^2-1/3$ is even and has zero integral against $U_1$. Moreover, it integrates to zero in $x$ against $U_{1/4}$. Hence these kernels define a martingale coupling of $(U_{1/4},U_1)$, whose conditional second moment is
\[
\frac{1}{3}-x\Big(\int y^4\,U_1(dy)-\frac{1}{3}\int y^2\,U_1(dy)\Big)=\frac{1}{3}-\frac{4}{45}x.
\]
In contrast with the martingale couplings of $(\mu_n,\nu_n)$, the conditional second moment now depends on $x$. It is smaller where the weight $1+x$ in the cost is larger, which lowers the cost. Indeed, using $\int x\,U_{1/4}(dx)=0$ and $\int x^2\,U_{1/4}(dx)=1/48$, the cost of this coupling is
\[
\int(1+x)\Big(\frac{1}{3}-\frac{4}{45}x-x^2\Big)\,U_{1/4}(dx)=V_Q-\frac{4}{45}\int x^2\,U_{1/4}(dx)=V_Q-\frac{1}{540}.
\]
Hence $V(\mu,\nu)\leq V_Q-1/540$, which together with $V(\mu_n,\nu_n)\to V_Q$ proves the claim.

\medskip

Every pair $(\mu_n,\nu_n)$ in Example~\ref{ex:parabola} satisfies~\eqref{ass:A} with the same $\kappa$, and is therefore irreducible, see Remark~\ref{rem:irreducible}. Moreover, no $\nu_n$ is carried by a hyperplane, and only the limit $\nu$ is. So neither~\eqref{ass:A} nor irreducibility alone gives the continuity of $V$ in dimension $d\geq2$, and the non-degeneracy of the limit $\nu$ in Theorem~\ref{thm:main} cannot be dropped.

\begin{remark}[Relation with the construction of \cite{BJ22}]\label{rem:BJ}
Example~\ref{ex:parabola} follows the mechanism of \cite{BJ22}. There, for a fixed integer $m \geq 2,$ $\mu_m$ is the uniform measure on the $m$ points $(i,0)$, $i=1,\dots,m$, and $\nu_{m,n}$ is obtained by splitting every atom of $\mu_m$ into two atoms of equal mass, at distance one along a line which makes the angle $\pi/(2n)$ with the horizontal axis. The projection along this direction onto the horizontal axis maps $\mu_m$ and $\nu_{m,n}$ to the same law. By strict convexity, this forces every martingale coupling to transport mass only along lines with this direction, so that the martingale coupling is unique, see \cite[Lemma~1.1]{BJ22}. As $n\to\infty$, these lines become horizontal, the constraint disappears, and the limit pair admits martingale couplings with smaller cost, see \cite[Proposition~1.2]{BJ22}. In both constructions, the target is concentrated near a set which flattens onto the horizontal axis, and the martingale condition in the vertical direction creates a constraint which is lost in the limit. In particular, also the limit target of \cite{BJ22} is carried by a line. The difference lies in the curvature. In \cite{BJ22} the vertical coordinate of the target is an affine function of the horizontal one along each line, and the constraint fixes the direction of transport. As a consequence, every point sends its mass to only two atoms of the target, and~\eqref{ass:A} fails for every $\kappa$. In Example~\ref{ex:parabola}, the target lies in a thin band
around a quadratic graph, and the martingale constraint
asymptotically fixes the conditional second moment of every
kernel. This leaves enough freedom to construct a martingale
coupling with density bounded away from zero and infinity,
so that~\eqref{ass:A} holds with a constant $\kappa$
independent of $n$.
\end{remark}

\section*{Acknowledgments}
F.M.~thanks Martin Huesmann for his kind hospitality during a visit to the Universität Münster and for fruitful discussions on martingale optimal transport, especially of the counterexamples discussed in the introduction.
F.M.~acknowledges support from the ERC Advanced Grant NEITALG, grant agreement No.~101198055. 
J.W.~acknowledges support from the Novo Nordisk and Villum Foundation.

AI tools were used for proofreading and for literature review. The authors take full responsibility for the final manuscript.

\begin{center}
  \FundingLogos
  
  \vspace{0.5em}
  \begin{tcolorbox}\centering\small
   
    Funded by the European Union. Views and opinions expressed are, however, those of the author(s) only and do not necessarily reflect those of the European Union or the European Research Council Executive Agency. Neither the European Union nor the granting authority can be held responsible for them. This project has received funding from the European Research Council (ERC) under the European Union’s Horizon Europe research and innovation programme (grant agreement No.~101198055, project acronym NEITALG).
    
  \end{tcolorbox}
\end{center}

\bibliographystyle{abbrv}
\bibliography{motpaper}

\end{document}